\documentclass[namedate,webpdf,imanum]{ima-authoring-template}%

\graphicspath{{Fig/}}

\theoremstyle{thmstyletwo}%
\newtheorem{thm}{Theorem}
\newtheorem{prop}[thm]{Proposition}%
\newtheorem{lem}[thm]{Lemma}
\newtheorem{rem}{Remark}%

\numberwithin{equation}{section}

\newcommand\bfd{{\mathbf d}}
\newcommand\bfe{{\mathbf e}}

\newcommand\bfr{{\mathbf r}}

\newcommand\bfu{{\mathbf u}}

\newcommand\bfw{{\mathbf w}}

\newcommand\bfA{{\mathbf A}}

\newcommand\bfE{{\mathbf E}}

\newcommand\bfK{{\mathbf K}}
\newcommand\bfM{{\mathbf M}}
\newcommand\bfW{{\mathbf W}}

\newcommand\bftheta{{\boldsymbol \vartheta}}

\def            \d          {\hspace{1.5pt}\mathrm{d}}
\newcommand{\R}{\mathbb{R}}
\DeclareMathOperator{\Span}{span}
\DeclareMathOperator{\Real}{Re}

\newcommand\eps{{\epsilon}}

\newcommand{\TbfM}{\leftidx{^\top}{\mathbf{M}}{}}
\newcommand{\TbfA}{\leftidx{^\top}{\bfA}{}}
\newcommand{\TbfK}{\leftidx{^\top}{\bfK}{}}
\newcommand{\bfde}{\mathbf{\dot e}}

\newcommand{\andquad}{\quad \textnormal{ and } \quad}

\newcommand{\Ga}{\varGamma}
\newcommand{\Om}{\varOmega}
\renewcommand{\nu}{\textnormal{n}}

\usepackage{tikz}
\usetikzlibrary{arrows}
\usepackage{adjustbox}

\usepackage{leftidx}

\usepackage{pgfplots}

\begin{document}

\DOI{DOI HERE}
\copyrightyear{2024}
\vol{00}
\pubyear{2024}
\access{Advance Access Publication Date: Day Month Year}
\appnotes{Paper}
\copyrightstatement{Published by Oxford University Press on behalf of the Institute of Mathematics and its Applications. All rights reserved.}
\firstpage{1}


\title[Error estimates for full discretization of Cahn--Hilliard equation with dyn.~b.c.]{Error estimates for full discretization of Cahn--Hilliard equation with dynamic boundary conditions}

\author{Nils Bullerjahn*
}
\author{Balázs Kovács\ORCID{0000-0001-9872-3474}
\address{\orgdiv{Institute of Mathematics}, \orgname{Paderborn University}, \orgaddress{\street{Warburger Str. 100}, \postcode{33098}, \state{Paderborn}, \country{Germany}}}}

\authormark{N.~Bullerjahn and B.~Kov\'acs}

\corresp[*]{Corresponding author: \href{email:bullerja@math.uni-paderborn.de}{bullerja@math.uni-paderborn.de}}

\received{29. }{7}{2024}
\revised{Date}{0}{Year}
\accepted{Date}{0}{Year}


\abstract{A proof of optimal-order error estimates is given for the full discretization of the Cahn--Hilliard equation with Cahn--Hilliard-type dynamic boundary conditions in a smooth domain. The numerical method combines a linear bulk--surface finite element discretization in space and linearly implicit backward difference formulae of order 1 to 5 in time. Optimal-order error estimates are proven. The error estimates are based on a consistency and stability analysis in an abstract framework, based on energy estimates exploiting the anti-symmetric structure of the second-order system.}


\maketitle

\section{Introduction}
In this paper we study the full discretization of the Cahn--Hilliard equation with Cahn--Hilliard-type dynamic boundary conditions, we discretize by linear bulk--surface finite elements in space and linearly implicit backward difference formulae of order 1 to 5 in time. We prove optimal-order error estimates between the numerical approximation and a sufficiently regular solution, see Theorem~\ref{Thm:ConvergenceMainResult} and the following discussion, for the combined bulk and surface $L^2$- and $H^1$-norms for the phase-field and chemical potential.

The Cahn--Hilliard equation was developed as a model for the continuous phase separation process of binary mixtures in a given domain in \cite{CahnHilliard1958}. The classical choice of boundary conditions are of Neumann type, which ensures the mass conservation, but force the free interface between the two components to intersect the boundary of the domain in a static angle of $\pi/2$, see, e.g., \cite{Jacqmin2000}. 

The short-ranged interactions between the boundary and the mixture are modeled by introducing an additional surface free energy which leads to different dynamic boundary conditions, depending on how one realizes the conservation of mass and energy dissipation: either separately for the bulk and surface, or in a combination. This leads to the following types of dynamic boundary conditions:
\begin{itemize}
	\item Allen--Cahn-type dynamic boundary conditions, \cite{FischerMaassDietrich1997}, satisfying conservation of the bulk mass and a dissipation of the total free energy.
	\item The GMS-model\footnote{The abbreviation represents the surnames of the respective authors.\label{footnote:abbreviation models}}, \cite{GoldsteinMiranvilleSchimperna2011}, satisfying conservation of the total mass and dissipation of the total free energy. The GMS-model is also referred as Cahn--Hilliard-type dynamic boundary conditions.
	\item The LW-model\footref{footnote:abbreviation models}, \cite{LiuWu2019}, satisfying conservation of the bulk and surface mass separately and dissipation of the total energy.
	\item Reaction rate dependent dynamic boundary conditions, \cite{KnopfLamLiuMetzger2021}, which interpolates the GMS-model and LW-model by a boundary coupling depending on a reaction rate.
\end{itemize}
\vspace{-2mm}
These types of dynamic boundary conditions all allow for a variable contact angle for the interface, which was also achieved in \cite{QianWangSheng2006}, via generalized Navier boundary conditions. 

A nice review paper collecting results on the Cahn--Hilliard equation with dynamic boundary conditions up to 2022 is \cite{Wu2022}.  

\smallskip
Works devoted to the discretization of the Cahn--Hilliard equation \emph{with} dynamic boundary conditions are: In \cite{CherfilsPetcuPierre2010} and \cite{CherfilsPetcu2014}, the authors have, respectively, considered Allen--Cahn and Cahn--Hilliard-type dynamic boundary conditions on a 2- or 3-dimensional slab. They have shown optimal-order fully discrete error estimates for a backward Euler time discretization. The recent work \cite{HarderKovacs2022} showed optimal-order semi-discrete error estimates for the GMS-model using a linear bulk--surface finite element approximation. The weak convergence up to a subsequence of the Cahn--Hilliard equation with reaction rate dependent dynamic boundary conditions for a bulk--surface finite element/backward Euler scheme was proved in \cite{KnopfLamLiuMetzger2021,Metzger2021}, and in \cite{Metzger2023} for an SAV scheme.  For the LW-model and a semi-discrete in time stabilized linearly implicit BDF2 scheme, an energy dissipation and error estimates are shown in \cite{MengBaoZhang2023}, and in \cite{LiuShenZheng2024} this scheme is used for various dynamic boundary conditions in a multiple scalar auxiliary variables (MSAV) approach. In \cite{AltmannZimmer2023} the authors derive a dissipation-preserving scheme by a PDAE formulation for Allen--Cahn-type dynamic boundary conditions, the LW-model, and GMS-model, and a semi-discrete in time backward Euler approach with convex-concave splitting. Additional structure-preserving schemes for the LW-model and GMS-model can be found in \cite{OkumuraFukao2024}.

%
%

In this paper we extend the results of \cite{HarderKovacs2022}. We study the full discretization of the weak formulation, based on a reformulation as a second-order system, using linear bulk--surface finite elements in space and $q$-step backward difference formulae (BDF method) in time, for $q=1,\dotsc,5$. We prove optimal-order fully discrete error estimates, which only require the potentials in the equation to satisfy a local Lipschitz condition for their first to third derivatives. We perform the stability and consistency analysis in an abstract functional analytical setting (similar to \cite{KovacsLubich2017,HarderKovacs2022}).
To our knowledge there are no fully discrete error estimates available for the Cahn--Hilliard equation with Cahn--Hilliard-type dynamic boundary conditions.

The fundamental technique of analyzing an isoparametric bulk--surface finite element method for an elliptic partial differential equation in a curved domain was introduced in \cite{ElliottRanner2013}, and many general results are proved and summarized in \cite{ElliottRanner2021}, even for a possibly evolving domain. General linear and semi-linear parabolic equations of second order with dynamic boundary conditions were analyzed in \cite{KovacsLubich2017}, providing optimal-order error estimates for a large class of problems in an abstract setting. 
The spatial discretization and its error estimates for our scheme are based on the ideas developed in these papers, and the consistency analysis will heavily depend on the general geometric approximations error estimates developed therein. 
Some early results for finite element discretizations of parabolic problems with dynamic boundary conditions on polyhedral domains were shown in \cite{Fairweather}.


\smallskip
A key contribution of this work is extending the basic idea for the semi-discrete stability analysis for the GMS-model \cite[Figure~1]{HarderKovacs2022}, for the Cahn--Hilliard equation on evolving surfaces \cite[Figure~1]{CHsurf}, and for Willmore flow \cite[Figure~5.2]{KovacsLiLubich2021} -- which are all based on the same core idea using energy estimates exploiting the anti-symmetric structure of the systems -- to a BDF full discretization (see Figure~\ref{fig:Scheme} herein).
Making a clear step towards the fully discrete error analysis for the arguably more involved Willmore flow.

In order to adapt the energy estimates for the $q$-step BDF method, we use a method combining results from the $G$-stability theory of \cite{Dahlquist1978} and the multiplier technique of \cite{NevanlinnaOdeh1981}. For PDEs, this approach -- testing with the errors -- was first used for linear parabolic differential equations on evolving surfaces in \cite{LubichMansourVenkataraman2013}, and later on for abstract quasi-linear parabolic problems in \cite{AkrivisLubich2015}. Fully discrete energy estimates -- testing with the discrete time derivatives of the errors -- were developed for mean curvature flow in \cite{KovacsLiLubich2019}. 
The above referenced energy estimates and, hence, the present work requires testing with the errors as well as their (discrete) time derivatives.

\smallskip
The paper is structured as follows:
In Section~\ref{section:Cahn--Hilliard intro} we introduce the Cahn--Hilliard equation and the GMS-model and its weak formulation as a system of second-order, which will be the basis of our discretization.
In Section~\ref{section:Semi-discretization CH} we formulate the semi-discretization in space by bulk--surface finite elements, in Section~\ref{section:BDF} the $q$-step backward difference method and in Section~\ref{section:MainResult} we state and discuss our main convergence result.
In Section~\ref{section:Stability} we establish the stability estimates, Section~\ref{section:Consistency} is devoted to the consistency analysis, and Section~\ref{section:MainProof} will combine these results to prove our main result.
In Section~\ref{section:NumericalExperiments} we present some numerical experiments to illustrate and complement our theoretical results.

\section{Cahn--Hilliard equation with dynamic boundary conditions\protect\footnote{Up to some minor changes, the text of this preparatory section is taken from \cite[Section 2]{HarderKovacs2022}.}}
\label{section:Cahn--Hilliard intro}

Let the bulk $\Om \subset \mathbb{R}^d$, $d=2,3$, be a bounded domain with an (at least) $C^2$ boundary $\Ga = \partial \Om$, which we call the surface. We denote by $\partial_\nu u$ the normal derivative of $u$ on $\Ga$, with $\nu$ being the unit outward pointing normal vector to $\Ga$. The surface gradient $\nabla_{\Ga} u$ on $\Ga$, of a function $u\colon \Ga \mapsto \mathbb{R}$, is given by $\nabla_{\Ga} u:= \nabla \bar u - (\nabla \bar u \cdot \nu )\nu$, where $\nu$ is the normal vector and $\bar u$ is an arbitrary extension of $u$. The Laplace--Beltrami operator on $\Ga$ is defined by $\varDelta_{\Ga} u = \nabla_\Ga \cdot \nabla_\Ga u$. For a more explicit definition of these operators, see, e.g., \cite[Section 2.2]{DziukElliott2013}. Let $\gamma u$ denote the trace of $u$ at $\Ga$ (see, e.g., \cite[Section 5.5]{Evans1998}), and let us denote the temporal derivative by ~$\dot{ } =\frac{\d}{\d t}$. For brevity, we write throughout the paper $\nabla_{\Ga}u$ instead of $\nabla_{\Ga}(\gamma u)$, and often $u$ instead of $\gamma u$.

In this paper we consider the Cahn--Hilliard equation subjected to the boundary conditions of the GMS-model, first derived in \cite{GoldsteinMiranvilleSchimperna2011}, which we will call dynamic boundary conditions of Cahn--Hilliard-type. A classical solution to this problem is a function $u\colon \Om \times [0,T] \to \mathbb{R}$ which solves the fourth-order equation:
\begin{subequations} \label{eq:strongCH}
	\begin{align}
		\dot u =&\ \varDelta \big(-\eps\varDelta u + \frac1\eps W'_{\Om}(u)\big) ~ &\text{ in } \Om,  \\
		\dot u =&\ \varDelta_\Ga \big(-\delta \kappa \varDelta_\Ga u +  \frac1\delta W'_{\Ga}(u)  + \eps \partial_\nu u\big) - \partial_\nu \big(-\delta \kappa \varDelta_\Ga u + \frac1\delta W'_{\Ga}(u) + \eps \partial_\nu u\big) ~ &\text{ on } \Ga,
	\end{align} 
\end{subequations}
with a sufficiently regular initial condition $u(0)=u^0$, and the diffusion interface parameters $\eps,\delta,\kappa>0$. The scalar functions $W_\Om$ and $W_\Ga$ are free energy potentials which are assumed to have locally Lipschitz first, second and third derivatives. A typical example is the double well potential $W(u)=(u^2-1)^2$, which satisfies these conditions.

\subsection{Weak formulation as a second-order system}
We introduce a chemical potential  $w\colon \bar \Om \times [0,T] \to \mathbb{R}$, to rewrite \eqref{eq:strongCH} into a system of second-order partial differential equations: Find $u,w \colon \bar \Om \times [0,T] \to \mathbb{R}$ satisfying
\begin{subequations} \label{eq:strong2CH}
	\begin{alignat}{2}
		\dot u =&\ \varDelta w &&\text{ in } \Om, \\
		w=&\ - \eps \varDelta u + \frac1\eps W'_\Om (u) &&\text{ in } \Om,
		\intertext{with Cahn--Hilliard-type dynamic boundary conditions }
		\dot u =&\ \varDelta_\Ga w - \partial_\nu w &&\text{ on } \Ga, \\
		w=&\ - \delta \kappa \varDelta_\Ga u + \frac1\delta W'_\Ga (u) + \eps \partial_\nu u \qquad &&\text{ on } \Ga.
	\end{alignat}
\end{subequations}

To define the weak formulation of problem \eqref{eq:strong2CH}, we use for $k\geq 0$ the standard Sobolev spaces $H^k(\Om)$ in the bulk (see, e.g., \cite[Section 5.2.2]{Evans1998}), and $H^k(\Ga)$ on the surface (see, e.g., \cite[Definition 2.11]{DziukElliott2013}). The variational formulation uses the Hilbert spaces 
\begin{equation}
\label{eq:Hilberst spaces and embedding}
	\begin{aligned}
		V=&\ \{ v \in H^1(\Om) ~|~ \gamma u \in H^1(\Ga)\} \andquad H=L^2(\Om) \times L^2(\Ga),\\
		&\ \text{with the dense embedding $V \to H$ : $\quad v \mapsto (v,\gamma v)$,}
	\end{aligned}
\end{equation}
and norms
\begin{equation}\label{eq:orignorms}
	\begin{aligned}
		\| u \|^2:=&\ \|u\|^2_V=\|u\|^2_{H^1(\Om)} + \|\gamma u\|^2_{H^1(\Ga)}, \\
		| (u_\Om , u_\Ga)|^2:=&\ \| (u_\Om , u_\Ga)\|^2_H = \|u_\Om\|^2_{L^2(\Om)} + \| u_\Ga \|^2_{L^2(\Ga)}. 
	\end{aligned}
\end{equation}
For an arbitrary $v\in V$ we often abbreviate pairs $(v,\gamma v) \in H$ by their first argument $v$. 

The weak formulation of the Cahn--Hilliard equation with dynamic boundary conditions of Cahn--Hilliard-type is obtained by testing \eqref{eq:strong2CH} with functions $\varphi^u , \varphi^w \in V$, integrating over the domain $\Om$ and using the classical Green's formula in $\Om$ to obtain
\begin{align*}
	\int_{\Om} \dot u \varphi^u =&\ - \int_\Om \nabla w \cdot \nabla \varphi^u + \int_\Ga \partial_\nu w \gamma \varphi^u, \nonumber \\
	\int_{\Om} w \varphi^w =&\  \eps \int_\Om \nabla u \cdot \nabla \varphi^w +\frac1\eps\int_\Om W'_{\Om}(u) \varphi^w - \eps\int_\Ga \partial_\nu u \gamma \varphi^w.
\end{align*}
Using the boundary conditions in \eqref{eq:strong2CH} in the boundary integrals and using Green's formula on the boundary (see, e.g., \cite[Theorem 2.14]{DziukElliott2013}) we obtain the weak problem: Find $u, w \in V$ such that for all $\varphi^u , \varphi^w \in V$ the following equations hold
\begin{subequations}
\label{eq:weakCHpara}
	\begin{alignat}{3}
		\bigg(\int_\Om \dot u \varphi^u + \int_\Ga \gamma \dot u \gamma \varphi^u\bigg) + \bigg(\int_\Om \nabla w \cdot  \nabla \varphi^u + & \int_\Ga  \nabla_{\Ga} w \nabla_{\Ga} \varphi^u\bigg) & = &\ 0, \\
		\bigg(\int_\Om w \varphi^w + \int_\Ga \gamma w \gamma \varphi^w\bigg) + \bigg(\eps \int_\Om \nabla u \cdot  \nabla \varphi^w +& \delta \kappa \int_\Ga  \nabla_{\Ga} u \nabla_{\Ga} \varphi^w\bigg) 
		& = &\ \frac1\eps\int_\Om W'_{\Om}(u) \varphi^w + \frac1\delta\int_\Ga W'_{\Ga}(u)\gamma \varphi^w.
	\end{alignat}
\end{subequations}

From now on we set the parameters that are involved in the problem \eqref{eq:weakCHpara} to $\eps=\delta=\kappa=1$ as their values only enter as additional constants in the error analysis we will carry out. The resulting problem then becomes: Find $u, w \in V$ such that for all $\varphi^u , \varphi^w \in V$ the following equations hold
\begin{subequations}
	\label{eq:weakCH}
	\begin{alignat}{3}
		\bigg(\int_\Om \dot u \varphi^u + \int_\Ga \gamma \dot u \gamma \varphi^u\bigg) + \bigg(\int_\Om \nabla w \cdot  \nabla \varphi^u + & \int_\Ga  \nabla_{\Ga} w \nabla_{\Ga} \varphi^u\bigg) & = &\ 0, \\
		\bigg(\int_\Om w \varphi^w + \int_\Ga \gamma w \gamma \varphi^w\bigg) + \bigg(\int_\Om \nabla u \cdot  \nabla \varphi^w +& \int_\Ga  \nabla_{\Ga} u \nabla_{\Ga} \varphi^w\bigg) 
		& = &\ \int_\Om W'_{\Om}(u) \varphi^w + \int_\Ga W'_{\Ga}(u)\gamma \varphi^w.
	\end{alignat}
\end{subequations}

\subsection{Abstract formulation}
It is possible to write the variational formulation in a more general abstract setting, using bilinear forms on the Hilbert spaces $V$ and $H$. This setting was introduced in \cite[Section 2.1.1]{KovacsLubich2017} for parabolic problems with dynamic boundary conditions, and used in the Cahn--Hilliard setting in \cite{HarderKovacs2022}. We have the Gelfand triple (see, e.g., \cite{Thomee2006}), with dense embeddings
\begin{align*}
	V \hookrightarrow H\simeq H^\ast \hookrightarrow V^\ast,
\end{align*} 
with the continuous and dense embedding $V \hookrightarrow H$ given by $v \mapsto (v,\gamma v)$, cf.~\eqref{eq:Hilberst spaces and embedding}. As a consequence, the duality $\langle \cdot , \cdot \rangle_V\colon V^\ast \times V \to \mathbb{R}$ coincides with $m(\cdot,\cdot)$ on $H\times V$.

We introduce the bilinear forms $a\colon V \times V \to \mathbb{R}$, $m\colon H \times H \to \mathbb{R}$ and $a^\ast \colon V \times V \to \mathbb{R}$ which we, respectively, define by setting
\begin{align*}
	a(u,v) = &\ \int_{\Om} \nabla u \cdot \nabla v + \int_\Ga \nabla_\Ga u \cdot \nabla_\Ga v, \\
	m((u_\Om, u_\Ga), (v_\Om,v_\Ga)) = &\ \int_{\Om} u_\Om v_\Om + \int_\Ga u_\Ga v_\Ga, \\
	\text{and} \qquad a^\ast(\cdot,\cdot) = &\ a(\cdot,\cdot)+m((\cdot,\gamma \cdot),(\cdot,\gamma \cdot )).
\end{align*}
We notice that the above defined norms \eqref{eq:orignorms} can also be expressed by
\begin{align*}
	\| u \|^2=a^\ast(u,u), \andquad |(u_\Om,u_\Ga)|^2=m((u_\Om, u_\Ga),(u_\Om, u_\Ga)).
\end{align*}
For elements of $V$ embedded into $H$, with the above notational convention, we even have $|u|^2=m(u,u)$. Furthermore, the bilinear form $a$ generates a semi-norm on $V$, denoted $\|u\|_a^2:=a(u,u)$. Lastly, for the non-linear terms we make the notational convention
\begin{align}
	m(W'(u),\varphi^w) := \int_{\Om} W'_\Om(u) \varphi^w + \int_\Ga W'_\Ga(u) \gamma \varphi^w. \label{eq:NotConvW}
\end{align}

Using this functional analytic setting we rewrite the weak formulation \eqref{eq:weakCH} as follows: Find functions $u \in C^1([0,T],H) \cap L^2([0,T],V)$ and $w \in L^2([0,T],V)$ such that, for any time $0<t \leq T$, and all $\varphi^u , \varphi^w \in V$,
\begin{subequations} \label{eq:absCH}
	\begin{align}
		m(\dot u(t),\varphi^u) + a(w(t),\varphi^u)=&\ 0, \label{eq:absCH1}\\
		m(w(t),\varphi^w) - a(u(t),\varphi^w)=&\ m(W'(u(t)),\varphi^w), \label{eq:absCH2}
	\end{align}
\end{subequations}
for given initial data $u(0)=u^0 \in V$. 

We define the free energy of the system \eqref{eq:absCH}, as the sum of the Ginzburg--Landau (bulk) free energy and of a Ginzburg--Landau surface free energy, see \cite[Section 2]{GoldsteinMiranvilleSchimperna2011}, by
\begin{align}
	\mathcal{E}(u):= \frac{1}{2} a(u,u) +  m(W(u),1). \label{eq:GinzLand}
\end{align}
Then we directly obtain, by testing \eqref{eq:absCH1} with $w$ and subtracting \eqref{eq:absCH2} tested with $\dot u$:
\begin{align*}
	0=a(u,\dot u) + m(W'(u),\dot u) + a(w,w) = \frac{\d}{\d t} \big( \mathcal{E}(u) \big) + \| w\|_a^2,
\end{align*}
that is, the dissipation of the bulk--surface free energy $\mathcal{E}$.

The following well-posedness and regularity results were proven in \cite[Theorem 3.2--3.6 and Remark 3.13]{GoldsteinMiranvilleSchimperna2011}: Let us assume that the non-linear potentials $W'_\Om$ and $W'_\Ga$ are locally Lipschitz continuous, vanish at zero and infinity, and their derivative is bounded from below, and we have an initial value $u^0 \in V$ of finite free energy $\mathcal{E}(u(0))<\infty$, which is smooth enough such that $w(\cdot,t)\in V$. Then, for some final time $T>0$, there exist functions
\begin{equation*}
		u \in L^\infty([0,T];V) \cap H^1([0,T];V), 
		\andquad
		w \in L^2([0,T];V) \cap L^\infty([0,T];V),
\end{equation*}
which satisfy \eqref{eq:absCH} almost everywhere in $(0,T)$. 

The following $H^2$-regularity result can be obtained by an adaptation of the arguments in \cite[Section 4.2-4.4]{ElliottRanner2015}: For $u^0 \in H^2(\Om)$ with $\gamma u^0 \in H^2(\Ga)$,
\begin{equation*}
	\begin{aligned}
		&u \in L^\infty([0,T];H^2(\Om)), ~ &&\text{ with } ~ &&\gamma u \in L^\infty([0,T];H^2(\Ga)),  \\
		&w \in L^2([0,T];H^2(\Om)), ~ &&\text{ with } ~ &&\gamma w \in L^2([0,T];H^2(\Ga)).  
	\end{aligned}
\end{equation*}
For more details see Appendix~\ref{appendix}.

\section{Semi-discretization of Cahn--Hilliard equation with dynamic boundary conditions\protect\footnote{Up to some changes, the text of this section is taken from \cite[Section 3]{HarderKovacs2022}.}}
\label{section:Semi-discretization CH}

For the spatial discretization of the above problem we use linear finite elements in the bulk and on the surface. This construction was developed by \cite[Sections 4,5]{ElliottRanner2013}, see also \cite[Sections 4--7]{ElliottRanner2021}, and by \cite[Section 4]{KovacsLubich2017}. We briefly recall the construction of the discrete domain and finite element space in the following.

\subsection{Bulk--surface finite elements}
\label{section:bulk--surface FEM}
The domain $\Om$ is approximated by a triangulation $\Om_h$ with maximal mesh width $h$, such that the boundary $\Ga_h:= \partial \Om_h$ of the corresponding polyhedral domain $\Om_h$ is an interpolation of $\Ga$, i.e. the vertices are on $\Ga$. Let $h$ be small enough that for every point $x \in \Ga_h$, there is a unique point $p\in \Ga$, such that $x-p$ is orthogonal to the tangent space $T_p\Ga$ of $\Ga$ at $p$, for more details on the construction see \cite[Section 4]{ElliottRanner2013}. We consider a quasi-uniform family of such triangulations $\Om_h$ of $\Om$, see, e.g., \cite[Definition 4.4.13]{BrennerScott2008}. 

The finite element space $V_h \nsubseteq H^1(\Om)$, for the discrete bulk $\Om_h$, is defined by
\begin{align}
	V_h:=\Span \{\phi_1,\dotsc,\phi_N\} \label{eq:DefVh},
\end{align}
where $\phi_j$ are continuous, piecewise linear nodal basis functions satisfying, for each node $(x_k)^N_{k=1}$,
\begin{align*}
	\phi_j(x_k)=\delta_{jk} \qquad \text{ for } j,k=1,\dotsc,N.
\end{align*}

Note that this definition transfers to the boundary: The restrictions of the basis functions to the boundary $\Ga_h$ form a surface finite element basis over the discrete boundary $\Ga_h$, see \cite[Section 5]{ElliottRanner2013}.

The discrete tangential gradient $\nabla_{\Ga_h}$ is piece-wisely defined, analogously to $\nabla_\Ga$, see also \cite[Section 4]{DziukElliott2013}. The discrete trace operator $\gamma_h v_h$ is defined by the restriction of the continuous functions $v_h \in V_h$ onto $\Ga_h$, with the same notational conventions as before, e.g.~$\nabla_{\Ga_h}v_h:=\nabla_{\Ga_h}(\gamma_h v_h)$.

In order to compare discrete functions in $V_h$ to functions in $V$, the lift operator $\cdot^\ell \colon V_h \to V$ is defined. The general construction is carried out in detail in \cite[Section 4]{ElliottRanner2013}, for functions on the boundary $v_h\colon \Ga_h \to \mathbb{R}$ the lift is defined as
\begin{align*}
	v^\ell_h\colon \Ga \to \mathbb{R}, \qquad \text{ with }~ v^\ell_h(p)=v_h(x) \quad \forall p \in \Ga,
\end{align*}
where $x \in \Ga_h$ is the unique point on $\Ga_h$ with $x-p\perp T_p\Ga$. This is then extended into the bulk. Then the lifted finite element space is given as $V^\ell_h:= \{v^\ell_h | v_h \in V_h\}$.

\subsection{Discrete norms and bilinear forms}
The discrete counterparts of the bilinear forms $a$, $m$ and $a^\ast$ are defined, 
for all $u_h,v_h \in V_h$, by
\begin{align*}
	a_h(u_h,v_h)=&\ \int_{\Om_h} \nabla u_h \cdot \nabla v_h + \int_{\Ga_h} \nabla_{\Ga_h} u_h \cdot \nabla_{\Ga_h} v_h, \\
	m_h(u_h,v_h)=&\ \int_{\Om_h} u_h v_h + \int_{\Ga_h} \gamma_h u_h \gamma_h v_h, \\
	\text{and} \qquad a^\ast_h(\cdot,\cdot)=&\ a_h(\cdot,\cdot)+m_h(\cdot,\cdot).
\end{align*}

The discrete norms and the semi-norm on $V_h$ are then defined, analogously to their continuous counterparts, by
\begin{equation}
\label{eq:discNorms}
	\|u_h\|_h^2:= a_h(u_h,u_h) + m_h(u_h,u_h), \quad
	|u_h|_h^2:= m_h(u_h,u_h), 
	\andquad \|u_h\|_{a_h}^2:= a_h(u_h,u_h).
\end{equation}

\subsection{A Ritz map}
We use a Ritz map to map the solution of the original problem into the discrete space, and also to define the initial data for the discrete problem. The Ritz map $\widetilde{R}_h\colon V \to V_h$ is defined, for arbitrary $u \in V$, see \cite[Section 4]{KovacsLubich2017}, by 
\begin{align}
\label{eq:DefRitzmap}
	a^\ast_h(\widetilde{R}_h u, \varphi_h)=a^\ast(u,\varphi_h^\ell), \qquad \text{for all} \qquad \varphi_h \in V_h. 
\end{align}
This map is well defined, due to the ellipticity of the bilinear form $a^\ast_h$, via the Riesz representation theorem. Note that the bilinear forms $a^\ast$ and $a^\ast_h$ contain boundary terms, therefore the Ritz map is clearly influenced by the boundary. We denote the (lifted) Ritz map by $R_hu:=(\widetilde{R}_h u)^\ell \in V^\ell_h$.

\subsection{The semi-discrete problem}

The semi-discretization of problem \eqref{eq:weakCH} reads: Find $u_h \in C^1([0,T],V_h)$ and $w_h \in L^2([0,T],V_h)$, such that, for time $0<t\leq T$ and for all $\varphi^u_h, \varphi^w_h \in V_h$,
\begin{subequations} \label{eq:sdcCHw}
	\begin{align}
		m_h(\dot u_h(t),\varphi^u_h)+a_h(w_h(t),\varphi^u_h)=&\ 0, \label{eq:sdcCHw1}\\
		m_h(w_h(t),\varphi^w_h)-a_h(u_h(t),\varphi^w_h)=&\ m_h(W'(u_h(t)),\varphi^w_h), \label{eq:sdcCHw2}
	\end{align}
	where the  initial data is defined to be $u_h(0)=u^0_h:=\widetilde{R}_h u^0 \in V_h$. Note that for $w$ no initial data is prescribed, however a unique $w_h(0)$ is obtained by solving the elliptic equation \eqref{eq:sdcCHw2}.
\end{subequations}

By standard theory on ordinary differential equations for short-time existence, well-posedness over long times is obtained by proving an energy dissipation see \cite[Proposition~3.1]{HarderKovacs2022}, in addition $H^2$-regularity was shown by adapting a technique by \cite[Theorem~4.3]{ElliottRanner2015}, and see in the Appendix herein. In \cite{HarderKovacs2022}, see Theorem~4.1 therein, optimal-order semi-discrete error estimates were shown.

\subsection{Matrix--vector formulation of the semi-discrete problem}
We rewrite \eqref{eq:sdcCHw2} into a matrix--vector formulation. We collect the nodal values of $u_h(\cdot,t)$ and $w_h(\cdot,t)$ into the vectors $\bfu(t) \in \mathbb{R}^N$ and  $\bfw(t) \in \mathbb{R}^N$, respectively. We define the mass and stiffness matrices
\begin{equation*}
	\bfM|_{ij}= m_h(\phi_j,\phi_i) 
	\andquad
	\bfA|_{ij}= a_h(\phi_j,\phi_i), \qquad i,j=1,\dotsc,N, 
\end{equation*}
where $\phi_j \in V_h$ are the set of basis functions of $V_h$ as in \eqref{eq:DefVh}. For the non-linear terms we collect the nodal values of $W_\Om'(u_h(\cdot,t))$ into $\bfW_\Om'(\bfu(t))$, and of $W_\Ga'(u_h(\cdot,t))$ into $\bfW_\Ga'(\bfu(t))$. We remember the notational convention \eqref{eq:NotConvW} and analogously use the notation $\bfM\bfW'(\bfu(t))$ for the combination. Then the equations \eqref{eq:sdcCHw} can be equivalently written as
\begin{subequations}
\label{eq:sMCH}
	\begin{align}
		\bfM \mathbf{\dot u} (t) + \bfA \bfw(t) = &\ 0 \label{eq:sMCH1},\\
		\bfM \bfw (t) - \bfA \bfu(t) = &\ \bfM \bfW'(\bfu(t)) . \label{eq:sMCH2}
	\end{align}
\end{subequations}	
This is the scheme which is the starting point of our stability analysis.

\section{Linearly implicit backward difference time discretization \label{section:BDF}}
Let $q=1,\dotsc,5$, $\tau >0$ be the time step size, and $q\tau \leq t_n := n\tau \leq T$ be a uniform partition of the time interval $[0,T]$. We assume the starting values $\bfu^0, \dotsc, \bfu^{q-1} \in \mathbb{R}^N$ to be given. The discretized time derivative and the extrapolation for the non-linear term are, respectively, given by
\begin{equation}
\label{eq:BDF diff and extrapolation}
	\mathbf{\dot u}^n= \frac{1}{\tau} \sum_{j=0}^{q} \delta_j \bfu^{n-j},\qquad \text{and} \qquad \mathbf{\widetilde u}^n:= \sum_{j=0}^{q-1}\gamma_j \bfu^{n-1-j}, \qquad \text{for $n \geq q$}.
\end{equation}
where the coefficients are given by the expressions
\begin{subequations}\label{eq:BdfCoeff}
	\begin{align}
		\label{eq:BDF generating function}
		\delta(\zeta)=&\sum_{j=0}^q \delta_j \zeta^j=\sum_{l=1}^q \frac{1}{l} (1-\zeta)^l, \\
		\gamma(\zeta)=&\sum_{j=0}^{q-1}\gamma_j \zeta^j=\frac{(1-(1-\zeta)^q)}{\zeta}.
	\end{align}
\end{subequations}
The $q$-step BDF methods are of order $q$, they are $A$-stable for $q=1,2$, $A(\alpha_q)$-stable for $q=3,\dotsc,6$ with $\alpha_3 =86.03^\circ$, $\alpha_4=73.35^\circ$, $\alpha_5=51.84^\circ$,  and $\alpha_6 = 17.84^\circ$, respectively, and unstable for $q \geq 7$, see, e.g.,  \cite[Section V.2]{HairerWanner1996}, \cite[Section 2.3]{AkrivisLubich2015}, while see \cite{AkrivisKatsoprinakis2020} for the exact $\alpha_q$ values. The starting values can be precomputed using either a lower order method with smaller time step size, or a Runge--Kutta method of suitable order.

Then the numerical scheme gives an approximation $\bfu^n$ of $u(t_n)$ and $\bfw^n$ of $w(t_n)$, for $q\tau \leq n \tau \leq T$, by the fully discrete system of linear equations
\begin{subequations} \label{eq:unmodifieddCH}
	\begin{align}
		\bfM \mathbf{\dot u}^n + \bfA \bfw^n=&\ 0, \label{eq:unmodifieddCH1}\\
		\bfM \bfw^n - \bfA \bfu^n =&\ \bfM \bfW'(\mathbf{\widetilde u}^n) .  \label{eq:unmodifieddCH2}
	\end{align}
\end{subequations}

\begin{rem}
	The linearly implicit numerical scheme \eqref{eq:unmodifieddCH} is well posed. In each time-step the following linear equation is solved:
	\begin{align*}
		\begin{bmatrix}
			\delta_{q} \bfM & \tau \bfA \\
			-\tau \bfA & \tau \bfM
		\end{bmatrix}
		\begin{bmatrix}
			\bfu^n\\
			\bfw^n
		\end{bmatrix} = \begin{bmatrix}
		\sum_{j=1}^q \delta_{j} \bfM \bfu^{n-j}\\
		\tau \bfM \bfW'(\mathbf{\widetilde u}^n)
		\end{bmatrix},
	\end{align*}
	which is uniquely solvable by the properties of the mass and stiffness matrices $\bfM$ and $\bfA$. 
\end{rem}

\begin{rem} \label{rem:MassEnergy}
	Conservation of the discrete combined mass $\bfM \bfu^n = \bfM \bfu^0$ is obtained by testing \eqref{eq:unmodifieddCH1} with $1 \in \R^N$, for initial values with constant discrete combined mass.
\end{rem}

\subsection{A modified fully discrete problem} \label{subsec:modifiedEq}
In the stability proof we will use energy estimates by testing the error equations \eqref{eq:CHNumericalScheme} with the discrete derivative of the errors, where the error is defined as the difference of the numerical solutions $\bfu^n$ and $\bfw^n$, and the Ritz map \eqref{eq:DefRitzmap} of the exact solutions $\bfu^n_\ast$ and $\bfw^n_\ast$ (see also \eqref{eq:RitzmapExactSol}). Therefore, it is natural that the error in the initial values is involved in the stability bound, even in the stronger $H^1$-norm. 

For the stability analysis the initial values $\bfu^i$ for $0\leq i \leq q-1$ are freely chosen, and therefore we can impose the condition $\bfu^i=\bfu_\ast^i$,
which makes the initial errors in $\bfu$ vanish. 

Similarly, we would like to have $\bfw^i=\bfw_\ast^i$, in order to avoid proving an $O(h^2)$ bound for the $H^1$-norm of these initial errors (which would not hold for the interpolation or Ritz map). In contrast to the initial values in $\bfu$, the initial values in $\bfw$ are not freely chosen, but are obtained by solving the second equation \eqref{eq:unmodifieddCH2}. For the analogous discussion in the semi-discrete setting we refer to \cite[Section~3.5]{HarderKovacs2022}, or \cite[Section~3.4]{KovacsLiLubich2021}

In order to make sense of the second equation in the numerical scheme \eqref{eq:unmodifieddCH2} for the initial time steps $i=0, \dotsc, q-1$, we extend the definition of the extrapolation by the exact value in these initial cases (see also \eqref{eq:defExtrapol}). Therefore, the initial values in $\bfw$ are given by:
\begin{align*}
	\bfM \bfw^i = \bfA \bfu^i + \bfM \bfW'(\bfu^i) , \qquad \text{for $i=0, \dotsc, q-1$.}
\end{align*}
Subtracting the same equation for the Ritz projection of the exact solution $\bfw_\ast$, we arrive at the identity
\begin{align*}
	\bfM (\bfw^i-\bfw_\ast^i)= -\bfM \bfd_w^i,
\end{align*}
where, for $i=0, \dotsc, q-1$, we define $\bfd_w^i$ as the nodal values of the semi-discrete defect $d_w(t_i)$ at time $t_i=i\tau$, see \cite[equation~(6.5)--(6.6)]{HarderKovacs2022}. Note again that the semi-discrete defect does not fulfill optimal-order estimates in the $H^1$-norm. Therefore, we shift the second equation by $\bftheta^n$, to set $\bfw^i=\bfw_\ast^i$.

The modified algorithm gives an approximation $\bfu^n$ of $u(t_n)$ and $\bfw^n$ of $w(t_n)$, for $q\tau \leq n \tau \leq T$, by the fully discrete system of linear equations
\begin{subequations}
\label{eq:dCH}
	\begin{align}
		\bfM \mathbf{\dot u}^n + \bfA \bfw^n=&\ 0, \label{eq:dCH1}\\
		\bfM \bfw^n - \bfA \bfu^n =&\ \bfM \bfW'(\mathbf{\widetilde u}^n) + \bfM \bftheta^n,  \label{eq:dCH2}
	\end{align}
\end{subequations}
where we define 
\begin{equation}
\label{eq:defShift}
	\begin{alignedat}{3}
		\bftheta^i = &\ \bfd_w^i , & \qquad & \text{for} \quad i = 0, 1, \dotsc, q-1 , \\
		\bftheta^n = &\ \bfd_w^{q-1} , & \qquad & \text{for} \quad n = q, q+1, \dotsc ,
	\end{alignedat}
\end{equation}
that is, the shift is constant except for the initial time steps for $i = 0, 1, \dotsc, q-1$.

Then, for all $i = 0, 1, \dotsc, q-1$, since the $\bfu$ initial values satisfy $\bfu^i=\bfu_\ast^i$, by \eqref{eq:dCH2} we have that the initial values for $\bfw$ satisfy:
\begin{align*}
	\bfM \bfw^i=&\ \bfA \bfu^i + \bfM \bfW'(\bfu^i) + \bfM \bftheta^i = \bfA \bfu_\ast^i + \bfM \bfW'(\bfu_\ast^i) + \bfM \bfd_w^i = \bfM \bfw_\ast^i, 
\end{align*}
that is -- as desired -- we have
\begin{equation*}
	\bfw^i = \bfw_\ast^i , \qquad \text{for $i = 0, 1, \dotsc, q-1$.}
\end{equation*}

\section{Main results: optimal-order fully discrete error estimates}
\label{section:MainResult}

The main goal of this paper is to establish optimal-order error estimates in space and time for the finite element discretization in space and linearly implicit $q$-step backward difference method in time for the Cahn--Hilliard equation with Cahn--Hilliard-type dynamic boundary conditions. In particular, to extend the semi-discrete error bounds of \cite{HarderKovacs2022} to the full discretization.


\begin{thm} \label{Thm:ConvergenceMainResult}
	Let $u$ and $w$ be sufficiently smooth solutions (see \eqref{eq:RegularityAss}) of the Cahn--Hilliard equation with Cahn--Hilliard-type dynamic boundary conditions \eqref{eq:strong2CH}, with nonlinear potentials satisfying \eqref{RegPotential}. Then, there exists $h_0 > 0$ such that for all $h \leq h_0$ and $\tau > 0$, satisfying the mild step size restriction $\tau^q \leq C_0 h^{2}$ (where $C_0 > 0$ can be chosen arbitrarily), the error between the solutions $u$ and $w$ and the fully discrete solutions $u_h^n$ and $w_h^n$ of \eqref{eq:dCH}, using linear bulk--surface finite element discretizations and linearly implicit backward difference time discretization of order $q=1,\dotsc,5$, satisfy the optimal-order error estimates, for $q\tau \leq t_n=n\tau \leq T$:
	\begin{equation*}
		\begin{aligned}
			&\ \| (u^n_h)^\ell - u(\cdot, t_n) \|_{L^2(\Om)} + \| \gamma((u^n_h)^\ell - u(\cdot, t_n)) \|_{L^2(\Ga)} \\
			&\ + h \big( \| (u^n_h)^\ell - u(\cdot, t_n) \|_{H^1(\Om)} + \| \gamma((u^n_h)^\ell - u(\cdot, t_n)) \|_{H^1(\Ga)}\big) \leq C \, (h^2+\tau^q), \\
			\text{and} \qquad 
			&\ \| (w^n_h)^\ell - w(\cdot, t_n) \|_{L^2(\Om)} + \| \gamma((w^n_h)^\ell - w(\cdot, t_n)) \|_{L^2(\Ga)} \\
			&\ + h \big( \| (w^n_h)^\ell - w(\cdot, t_n) \|_{H^1(\Om)} + \| \gamma((w^n_h)^\ell - w(\cdot, t_n)) \|_{H^1(\Ga)}\big) \leq C \, (h^2+\tau^q). 
		\end{aligned}
	\end{equation*}
	Further, the errors in the discrete time derivative \eqref{eq:BDF diff and extrapolation}, for finite element functions denoted here by $\partial^\tau_q$, in $u$ satisfy, for $q\tau \leq t_n=n\tau \leq T$:
	\begin{equation*}
		\begin{aligned}
			\bigg( \tau \sum_{k=q}^n \Big(&\| (\partial^\tau_q u^k_h)^\ell - \dot u(\cdot,t_k) \|_{L^2(\Om)}^2 + \| (\partial^\tau_q u^k_h)^\ell - \dot u(\cdot,t_k) \|_{L^2(\Ga)}^2 \\
			& + h \, \big(\| (\partial^\tau_q u^k_h)^\ell - \dot u(\cdot,t_k) \|_{H^1(\Om)}^2 + \| (\partial^\tau_q u^k_h)^\ell - \dot u(\cdot,t_k) \|_{H^1(\Ga)}^2\big) \Big) \bigg)^{\frac{1}{2}} \leq C \, (h^2+\tau^q).
		\end{aligned}
	\end{equation*}
	The constant $C>0$ depends on the Sobolev norms of the solutions in accordance with \eqref{eq:RegularityAss}, and exponentially on the final time $T$, but is independent of $h$ and $\tau$. 
\end{thm}

\begin{rem}
	(i) The $h_0$ depends on geometric approximation results (see, e.g., Section~\ref{section:bulk--surface FEM} and Section~\ref{section:Consistency}). The requirement on the smallness of both parameters is merely technical. 
	
	(ii) The mild step size restriction is needed to establish an $L^\infty$-bound on the error in the stability analysis. This lets us reduce the assumptions on the nonlinear potentials $W_\Om$ and $W_\Ga$ to a \emph{local} Lipschitz condition, which covers important examples such as the double well potential. The restriction can be even weakened to $\tau^q \leq C_0 h^{3/2+\eps_0}$, where $\eps_0,C_0 > 0$ can be chosen arbitrarily, as is apparent from the assumptions in Proposition~\ref{Prop:StabilityEst}.
	
	(iii) The constant $C$ also depends on a polynomial expression of the inverses of $\eps$ and $\delta$, that is, limit processes for the interface parameters are not covered by our error estimates.
\end{rem}

For Theorem~\ref{Thm:ConvergenceMainResult} sufficient regularity assumptions on the nonlinear potentials are,
\begin{equation}
\label{RegPotential}
	W_\Om^{(k)} \text{ and } W_\Ga^{(k)} \text{ are locally Lipschitz continuous for }k=0,\cdots,3, 
\end{equation}
and sufficient regularity assumptions on the exact solutions are: 
\begin{equation} \label{eq:RegularityAss}
	\begin{alignedat}{3}
		&\ u \in  H^{q+3}([0,T];H^2(\Om)), & \quad &
		\gamma u \in  H^{q+3}([0,T];H^2(\Ga)), \\
		&\ w \in  H^{3}([0,T];H^2(\Om)), & \quad &
		\gamma w \in  H^{3}([0,T];H^2(\Ga)).
	\end{alignedat}
\end{equation}

Recall that by standard theory we have, for any $u \in H^1([0,T];X)$, that $u \in C([0,T],X)$ such that
\begin{align*} 
	\max_{0\leq t \leq T} \| u(t)\|_X \leq c(T) \|u\|_{H^1([0,T];X)},
\end{align*}
see, e.g., \cite[Section 5.9.2]{Evans1998}. Since the Sobolev embedding $H^2 \subset L^\infty$ is continuous in dimensions $d=2,3$, see, e.g., \cite[Section 5.6.3]{Evans1998}, we have the regularity 
\begin{equation*}
	\begin{aligned}
		u &\in C^1([0,T];H^1(\Om)) \cap W^{q+2,\infty}([0,T];L^{\infty}(\Om)), \quad &&\gamma u \in C^1([0,T];H^1(\Ga)) \cap W^{q+2,\infty}([0,T];L^{\infty}(\Ga)), \\
		w &\in C([0,T];H^1(\Om)) \cap W^{2,\infty}([0,T];L^\infty(\Om)), \quad &&\gamma w \in  C([0,T];H^1(\Ga)) \cap W^{2,\infty}([0,T];L^\infty(\Ga)),
	\end{aligned}
\end{equation*}
where the first spaces always indicate the regularity required for the weak formulation, while the second spaces are required for the error analysis.

In the following sections we will prove this result by a combination of a stability and consistency analysis. 

In order to achieve the stability bound, we follow the general scheme of the stability proof in \cite[Section 5]{HarderKovacs2022}, first developed in \cite{KovacsLiLubich2021} for Willmore flow, which has a similar anti-symmetric structure. We derive two sets of energy estimates by exploiting the anti-symmetric structure of \eqref{eq:dCH}: (i) In the first step an energy estimate for the error in $u$ is proved, which comes with a critical term involving the discrete time derivative of the error in $u$. (ii) For the second energy estimate we take the discrete time derivative of the second equation, which leads to a bound for the critical term in (i) and a bound for the error in $w$. The fully discrete energy estimates rely on the $G$-stability theory of \cite[Theorem 3.3]{Dahlquist1978} and the multiplier technique of \cite[Section 2]{NevanlinnaOdeh1981}. This method has already proved to be effective for different partial differential equations, see, e.g., \cite{LubichMansourVenkataraman2013},  \cite{AkrivisLubich2015}, \cite{KovacsLiLubich2019}, \cite{LLG}, \cite{HMHF}, \cite{ContriKovacsMassing2023}. 

In the consistency part we are concerned with the approximation error which comes from the discretization and enters as defect terms in the stability analysis. 
The consistency analysis relies on an interplay of the error estimates for the backward time differentiation and the error estimates for the spatial discretization \cite[Section~6]{HarderKovacs2022} (which rely on error estimates of the nodal interpolations \cite[Proposition~5.4]{ElliottRanner2013}, error estimates for a Ritz map \cite[Lemma~3.1,3.3]{KovacsLubich2017}, and geometric approximation errors in the bilinear forms \cite[Lemma 4.2]{KovacsLubich2017}). 
To obtain the optimal-order fully discrete consistency bounds, the major difficulty stems from the fact that in the stability estimate the discrete derivative of the fully discrete defect terms appear. In order to deal with these terms, we have to extend our notion of the discrete derivative by a sufficiently regular extrapolation to the whole time interval. This is done by a quintic Hermite interpolation. 

\section{Stability}
\label{section:Stability}

\subsection{Error equations for the full discretization}
Let us denote the Ritz map \eqref{eq:DefRitzmap} of the exact solutions, at any time $t_n=n\tau \in [0,T]$, by
\begin{align}
	u^n_\ast=\widetilde{R}_h u(\cdot,t_n)\qquad \text{and} \qquad w^n_\ast=\widetilde{R}_h w(\cdot,t_n). \label{eq:RitzmapExactSol}
\end{align}

We collect the nodal values into $\bfu_\ast^n$ and $\bfw_\ast^n$, then these vectors satisfy the numerical scheme for $n\geq q$, only up to some defects $\bfd^n_u , \bfd^n_w$:
\begin{subequations} \label{eq:exdCH}
	\begin{align}
		\bfM \mathbf{\dot u}^n_\ast + \bfA \bfw^n_\ast=&\  \bfM \bfd^n_u, \\
		\bfM \bfw^n_\ast - \bfA \bfu^n_\ast =&\ \bfM \bfW'(\mathbf{\widetilde u}^n_\ast) + \bfM \bfd^n_w .
	\end{align}
\end{subequations}

The errors between the discrete solution and the Ritz map of the exact solutions are denoted by $\bfe^n_u=\bfu^n-\bfu^n_\ast$ and $\bfe^n_w=\bfw^n-\bfw^n_\ast$. By subtracting the equations \eqref{eq:dCH} from \eqref{eq:exdCH} we obtain, that the errors satisfy the following error equations
\begin{subequations}
	\label{eq:CHNumericalScheme}
	\begin{align}
		\bfM \bfde^n_u + \bfA \bfe^n_w=&\  -\bfM\bfd^n_u, \label{eq:CHNumericalScheme1}\\
		\bfM \bfe^n_w - \bfA \bfe^n_u =&\  \bfM \bfr^n, \label{eq:CHNumericalScheme2}
	\end{align}
	where, recalling \eqref{eq:defShift}, we set
	\begin{align}
		\bfr^n=\bfW'(\mathbf{\widetilde u}^n) - \bfW'(\mathbf{\widetilde u}^n_\ast) - \bfd^n_w + \bftheta^n. \label{eq:defr}
	\end{align}
\end{subequations}

The defects are estimated using a discrete dual norm of the space $V_h$, defined by
\begin{align*}
	\| d_h \|_{\ast,h} = \sup_{0\not= \varphi_h \in V_h} \frac{m_h(d_h,\varphi_h)}{\|\varphi_h\|_h},
\end{align*}
for any $d_h \in V_h$. We use the same notation for the dual norm for a nodal vector $\bfd$ corresponding to $d_h$.

For the initial time steps $t_i$, $i=0,\dotsc,q-1$, we define  the defects as the value of the defects in the semi-discrete problem \eqref{eq:sdcCHw}, i.e.~$\bfd^i_u$ and $\bfd^i_w$ are the vectors of nodal values of $d^u_h(i\tau)$ and $d^w_h(i\tau)$, see \cite[equation~(6.5)--(6.6)]{HarderKovacs2022}.

The stability result uses energy estimates, testing both with the errors and their time derivatives. For these estimates we use the matrix--vector formulation, and the following natural (semi-)norms for $\bfM$, $\bfA$, and $\bfK=\bfM+\bfA$, corresponding to the discrete norms in \eqref{eq:discNorms}, for any $\bfu \in \R^N$ corresponding to the finite element function $u_h \in V_h$,
\begin{align*}
	\|\bfu\|_\bfK^2 = &\ \bfu^T \bfK \bfu = \bfu^T \big( \bfM + \bfA \big) \bfu = \|u_h\|_h^2 , \\
	\|\bfu\|_\bfM^2 = &\ \bfu^T \bfM \bfu = |u_h|_h , \\
	\text{and} \qquad 
	\|\bfu\|_\bfA^2 = &\ \bfu^T \bfA \bfu = \|u_h\|_{a_h}^2 .
\end{align*}

\subsection{An $L^\infty$ bound for the Ritz map}
\label{subsec:Linftyu}

To control the locally Lipschitz continuous potential $\bfW'$ of the extrapolated value of the numerical and the exact solution, we show that all inputs are in a compact set, which allows us to use Lipschitz continuity for these terms. As a preparation we show here an $h$-uniform $L^\infty$-bound of the Ritz map of the exact solution $u_*^n = \widetilde{R}_h u(\cdot,t_n)$. 

By the inverse estimate (see, e.g., \cite[Theorem~4.5.11]{BrennerScott2008}), the $h$-uniform norm equivalence between the discrete and original $L^2$-norms (see Remark~\ref{rem:normequivalence}), the $L^2$-norm error estimates of the finite element interpolation operator $\widetilde I_h v \in V_h$, with lift $I_h v =(\widetilde I_h v)^\ell \in V_h^\ell$ (\cite[Proposition~2.7]{Demlow2009}, or Lemma~\ref{Lemma:InterpolationError} below), and the error estimates of the Ritz map (\cite[Lemma 4.4]{KovacsLubich2017}, or Lemma~\ref{Lemma:RitzMapError} below) we estimate, for any $0 \leq t \leq T$ (suppressed in the notation) and $j = 0, 1$ (letting $u^{(1)}:=\dot u$),
\begin{align*}
	\| \widetilde R_h u^{(j)} - \widetilde I_h u^{(j)} \|_{L^\infty(\Om_h)} \leq&\  c h^{-d/2} \| \widetilde R_h u^{(j)} - \widetilde I_h u^{(j)} \|_{L^2(\Om_h)} \\
	\leq&\  c h^{-d/2} \| R_h u^{(j)} - u^{(j)} \|_{L^2(\Om)} 
	+ c h^{-d/2} \| u^{(j)} - I_h u^{(j)} \|_{L^2(\Om)} \\
	\leq&\ c h^{2-d/2} \| u^{(j)} \|_{H^2(\Om)},
\end{align*}
and using the regularity of the exact solution \eqref{eq:RegularityAss}, and the $L^\infty$-stability of the linear finite element interpolation operator, for any $0 \leq t \leq T$ and $j = 0, 1$, we have
\begin{equation}
	\label{eq:exactuLinfty}
	\begin{aligned}
		\| u_\ast^{(j)} \|_{L^\infty(\Om_h)} 
		\leq &\ \| \widetilde R_h u^{(j)} - \widetilde I_h u^{(j)} \|_{L^\infty(\Om_h)} 
		+ c \| I_h u^{(j)} \|_{L^\infty(\Om)} 
		\\
		\leq&\  c h^{2-d/2} \| u^{(j)} \|_{H^2(\Om)} + c \| u^{(j)} \|_{L^\infty(\Om)} \\
		\leq&\ c.
	\end{aligned}
\end{equation}
We note here, that for higher order finite elements one should estimate $\| I_h u^{(j)} \|_{L^\infty(\Om)} \leq \| I_h u^{(j)} - u^{(j)} \|_{L^\infty(\Om)}  + \| u^{(j)} \|_{L^\infty(\Om)}$ and use an estimate on the interpolation operator, see e.g. \cite[Theorem~4.28]{ElliottRanner2021}.

For the extrapolation we simply have
\begin{align*}
	\| \mathbf{\widetilde u}^k_\ast \|_{L^\infty(\Om_h)} \leq \sum_{j=0}^{q-1} |\gamma_j|  \|\bfu_\ast^{k-1-j}\|_{L^\infty(\Om_h)} \leq c.
\end{align*}

In the stability proof we then argue by a bootstrapping argument that, since the extrapolation only depends on past values, $\mathbf{\widetilde u}$ is $L^\infty$-close to $\mathbf{\widetilde u}_\ast$, and hence all inputs are in a compact set.

\subsection{Results by Dahlquist and Nevanlinna \& Odeh}

In order to derive energy estimates for BDF methods up to order $5$, we need the following important results from the $G$-stability theory of \cite[Theorem~3.3]{Dahlquist1978} and the multiplier technique of \cite[Section 2]{NevanlinnaOdeh1981}.
\begin{lem}[{\cite[Theorem 3.3]{Dahlquist1978}}]
	\label{Lemma:DahlquistOrig} 
	Let $\delta (\zeta)= \sum_{j=0}^q \delta_j \zeta^j$ and $\mu(\zeta)=\sum_{j=0}^q \mu_j \zeta^j$ be polynomials of degree at most $q$, and at least one of them of degree $q$, that have no common divisor. Let $\langle\cdot,\cdot\rangle$ denote an inner product on $\mathbb{R}^N$ with associated norm $|\cdot|$. If
	\begin{align*}
		\Real \frac{\delta(\zeta)}{\mu(\zeta)} >0 , \qquad \text{ for all } \quad \zeta \in \mathbb{C}, \ |\zeta|<1,
	\end{align*}
	then there exists a symmetric positive definite matrix $G=(g_{ij}) \in \mathbb{R}^{q\times q}$ and real numbers $\gamma_0,\dotsc,\gamma_q$ such that, for all $w_0,\dotsc,w_q \in \mathbb{R}^N$, the following identity holds
	\begin{align*}
		\Real \Big\langle \sum_{i=0}^q \delta_i w_{q-i}, \sum_{j=0}^q \mu_j w_{q-j} \Big\rangle = \sum_{i,j=1}^q g_{ij} \langle w_i,w_j \rangle - \sum_{i,j=1}^q g_{ij} \langle w_{i-1},w_{j-1} \rangle +\Big|\sum_{i=0}^q \gamma_i w_i \Big|.
	\end{align*}
\end{lem}

Additionally, we will need this result for an arbitrary semi-inner product, therefore we use the following extension of the above result to semi-inner products by \cite[Lemma~3.5]{ContriKovacsMassing2023}.
\begin{lem} 
\label{Lemma:Dahlquist} 
	Let the polynomials $\delta(\zeta)$ and $\mu(\zeta)$ satisfy the conditions of Lemma~\ref{Lemma:DahlquistOrig}, and let $(\cdot , \cdot )$ be a semi-inner product on a Hilbert space $H$. Then, for the matrix $G=(g_{ij}) \in \mathbb{R}^{q\times q}$ from Lemma~\ref{Lemma:DahlquistOrig}, and for all $w_0,\dotsc,w_q \in H$ we have
	\begin{align*}
		 \Real \Big( \sum_{i=0}^q \delta_i w_{q-i}, \sum_{j=0}^q \mu_j w_{q-j} \Big) \geq \sum_{i,j=1}^q g_{ij} ( w_i,w_j ) - \sum_{i,j=1}^q g_{ij} ( w_{i-1},w_{j-1} ) .
	\end{align*}
\end{lem}

The application of $G$-stability to the BDF schemes, for $q=1,\dotsc,5$, is ensured by the following result.
\begin{lem}[{\cite[Section 2]{NevanlinnaOdeh1981}}]
	\label{Lemma:MultiplierTechnique}
	For $1\leq q\leq 5$, there exists $0 \leq \eta <1$ such that for $\delta (\zeta)= \sum_{j=0}^q \frac{1}{j} (1-\zeta)^j$,
	\begin{align*}
		\Real \frac{\delta(\zeta)}{1-\eta \zeta} >0 , \qquad \text{ for all } \quad \zeta \in \mathbb{C}, \ |\zeta|<1 .
	\end{align*}
	The classical values of $\eta$ from \cite[Table]{NevanlinnaOdeh1981} are found to be $\eta=0$, $0$, $0.0836$, $0.2878$, $0.8160$ for $q=1,\dotsc,5$ respectively.
\end{lem}
The exact multipliers were computed in \cite{AkrivisKatsoprinakis2015}, while multipliers for BDF-6 were derived in \cite{AkrivisChenYuZhou2021}.

We thus introduce the $G$-(semi-)norm associated to the (semi-)inner product $(\cdot,\cdot)$ on a Hilbert space $H$: Given a collection of vectors $W^n=(w^n,\dotsc,w^{n-q+1}) \subset H$, we define
\begin{align*}
	\|W^n\|_G^2:= \sum_{i,j=1}^q g_{ij} ( w^{n-i+1},w^{n-j+1} ),
\end{align*}
where $G$ is the symmetric positive definite matrix appearing in Lemma~\ref{Lemma:DahlquistOrig}  and \ref{Lemma:Dahlquist}. Then with the smallest and largest eigenvalues of $G$, denoted by $\lambda_0$ and $\lambda_1$, we have the inequalities:
\begin{align}
	\lambda_0 \|w^n\|^2 \leq \lambda_0 \sum_{j=1}^q \|w^{n-j+1}\|^2 \leq \|W^n\|_G^2 \leq \lambda_1 \sum_{j=1}^q \|w^{n-j+1}\|^2, \label{eq:GnormEst}
\end{align}
where $\|\cdot\|$ is the (semi-)norm on $H$ induced by the (semi-)inner product $(\cdot,\cdot)$.
Later on, an additional subscript, e.g.~$\|\cdot\|_{G,\bfM}$, specifies which (semi-)scalar product generates the $G$-weighted (semi-)norm.

These results have previously been applied to the error analysis for the BDF time discretization of PDEs when testing the error equation with the error in \cite{LubichMansourVenkataraman2013}, and the discrete time derivative of the error in \cite{KovacsLiLubich2019}.

We additionally use the following lemma, relating the BDF time derivative to the first-order backward difference, which stems from the zero-stability of the BDF method for $q\leq 6$, see \cite[Section III.3]{HairerNorsettWanner1993}.

\begin{lem} \label{Lemma:partialtodot}
	Let $(a^k) \in H$ be a sequence in a normed vector space (with norm $\| \cdot \|$). 
	Recalling the $q$-step BDF discrete derivative \eqref{eq:BDF diff and extrapolation}, and the first order backward difference formula:
	\begin{alignat*}{3}
		\partial^\tau_q a^k :=&\ \frac{1}{\tau} \sum_{j=0}^q \delta_j a^{k-j} , \qquad &&\text{for } &&k \geq q, \\
		\text{and} \qquad \partial^\tau a^k :=&\ \partial^\tau_1 a^k , \qquad &&\text{for } &&k \geq 1.
	\end{alignat*}
	Then, for a $\tau$- and $n$-independent $c > 0$, we have that
	\begin{align*}
		\sum_{k=q}^n \| \partial^\tau a^k \|^2 \leq  c \sum_{k=q}^n \| \partial^\tau_q a^k \|^2 + c \sum_{i=1}^{q-1} \| \partial^\tau a^i \|^2 .
	\end{align*}
\end{lem}
\begin{proof}
	This bound follows directly from a factorization of the generating polynomial $\delta(\zeta)= (1-\zeta) \sigma(\zeta)$ with $\sigma(\zeta)=\sum_{j=0}^{q-1} \sigma_j \zeta^j$, which allows us to write, for $q \leq k \leq n$,
	\begin{align}
		\partial^\tau_q a^k = \frac{1}{\tau} \sum_{j=0}^q \delta_j a^{k-j} = \sum_{j=0}^{q-1} \sigma_j \partial^\tau a^{k-j}. \label{eq:factorization}
	\end{align}	
	A detailed proof is given in \cite[equation (10.12)--(10.14)]{KovacsLiLubich2019}. \hfill
\end{proof}

\subsection{Stability bound}

We will now state and prove the fully discrete stability for the numerical method \eqref{eq:dCH}, which extends the underlying semi-discrete stability result \cite[Proposition~5.1]{HarderKovacs2022}.

\begin{prop}
\label{Prop:StabilityEst}
	Assume that the exact solution satisfies the regularity assumption \eqref{eq:RegularityAss}. Consider the linearly implicit BDF time discretization of order $1 \leq q \leq 5$. 
	Assume that, for step-sizes satisfying $\tau^q \leq C_0 h^\kappa$ (with arbitrary $C_0 > 0$), for a $\kappa > \frac{3}{2}$ the defects are bounded by 
	\begin{equation}
	\label{eq:assumed bounds defects}
		\begin{alignedat}{2}
			\|\bfd^k_u \|_{\ast,h} \leq&\ ch^\kappa , \qquad \quad\|\bfd^k_w \|_{\ast,h} &&\leq\ ch^\kappa, \qquad \| \bfd_w^{q-1} \|_{\ast,h}\leq ch^\kappa \\
			\|\mathbf{\dot d}^k_u \|_{\ast,h} \leq&\ ch^\kappa, \qquad \quad \|\mathbf{\dot d}^k_w \|_{\ast,h} &&\leq\ ch^\kappa, \qquad \| \partial^\tau \bfd_w^i \|_{\ast,h} \leq ch^\kappa ,
		\end{alignedat}
	\end{equation}
	for $q\tau \leq k \tau \leq T$ and $1 \leq i \leq q-1$. 
	Further assume that the errors of the starting values satisfy
	\begin{align}
	\label{eq:defI}
		I_h^{q-1} := \sum_{i=0}^{q-1} \Big( \|\bfe_u^i\|_\bfK^2 + \|\bfe_w^i\|_\bfK^2 \Big) + \tau \sum_{i=1}^{q-1} \| \partial^\tau \bfe_u^i \|_\bfK^2 + \tau \sum_{i=1}^{q-1} \| \partial^\tau \bfe_w^i \|_\bfK^2 \leq c h^{2\kappa} . 
	\end{align}
	Then there exist $h_0 > 0$, such that for all $h\leq h_0$ and with $\tau^q \leq C_0 h^\kappa$, the following stability estimate hold for \eqref{eq:CHNumericalScheme}, for all $q \tau \leq t_n=n \tau \leq T$,
	\begin{align}
		\| \bfe^n_u \|_\bfK^2 + \| \bfe^n_w \|_\bfK^2 + \tau \sum_{k=q}^n \| \bfe^k_w \|_\bfK^2 + \tau \sum_{k=q}^n \| \bfde^k_u \|_\bfK^2 \leq C\Big(I_h^{q-1} +  D_h^n \Big) , \label{eq:StabilityEst}
	\end{align}
	where the constant $C>0$ is independent of $h$, $\tau$ and $n$, but depends exponentially on the final time $T$, and where 
	the term $D_h^n$ on the right-hand side collects the defects:
	\begin{align}
		D_h^n:=&\ \| \bfd_u^n \|_{\ast,h}^2 + \| \bfd_u^q \|_{\ast,h}^2 + \| \bfd_w^{q-1} \|_{\ast,h}^2 +  \tau \sum_{k=q}^n \Big( \|\bfd^k_u \|_{\ast,h}^2 + \| \bfd^k_w \|_{\ast,h}^2 \Big) \nonumber \\
		&\ + \tau \sum_{k=q}^n \Big(\| \mathbf{\dot d}^k_u \|_{\ast,h}^2 + \| \mathbf{\dot d}^k_w \|_{\ast,h}^2\Big) + \tau \sum_{i=1}^{q-1}\|\partial^\tau \bfd_w^i \|_{\ast,h}^2 . \label{eq:defD}
	\end{align}
\end{prop}

\begin{proof}
The proof of the stability bound follows the general scheme of the stability proof in \cite[Section 5]{HarderKovacs2022}, which was first developed in \cite[Section 5]{KovacsLiLubich2021} for Willmore flow, and we extend this approach to the fully discrete case. 
The proof relies on energy estimates, obtained by testing the error equations with the fully discrete errors and their discrete time-derivatives, using Dahlquist's $G$-stability theory (Lemma~\ref{Lemma:DahlquistOrig} and \ref{Lemma:Dahlquist}) and the multiplier technique of Nevanlinna and Odeh (Lemma~\ref{Lemma:MultiplierTechnique}). The core idea of the proof is presented in the diagram in Figure~\ref{fig:Scheme}: (i) In the first step an energy estimate for $\bfe^n_u$ is proved, which comes with a critical term involving $\bfde^n_u$. (ii) In the second estimate we use the discrete BDF time derivative of \eqref{eq:CHNumericalScheme2}, which leads to a bound for this critical term and for $\bfe^n_w$ as well. 

\begin{figure}[hb]
	\begin{adjustbox}{width=\textwidth}
		\begin{tikzpicture}
			\node at (2.5,5) (erroreq1) {$\bfM \dot{\bfe}_u^n + \bfA \bfe_w^n = - \bfM \bfd_u^n$};
			\node at (2.35,4.5) (erroreq2) {$\bfM \bfe_w^n - \bfA \bfe_u^n = \bfM \bfr^n$};
			
			\node at (-0.95,-0.345) (energyest1) {\small $\displaystyle \|\bfe^n_u\|_\bfK^2 + \tau \! \sum \!\! ~\|\bfe_w^n\|_\bfK^2 \lesssim \! \underbrace{\tau \sum \!\! ~\|\dot{\bfe}_u^n\|_\bfK^2}_{\text{{\tiny critical term}}} +I_h^{q-1} + D_h^n$};
			\node at (8.75,-0.1) (energyest2) {\small $\displaystyle \|\bfe_w^n\|_\bfK^2 + \tau \! \sum \!\! ~\|\dot{\bfe}_u^n\|_\bfK^2 \lesssim I_h^{q-1}  +  D_h^n$};
			
			\pgfmathsetmacro{\yshftii}{-0.75}
			\node at (3,4.5+\yshftii) (dterroreq2) {$\bfM \dot{\bfe}_w^n - \bfA \dot{\bfe}_u^n = \bfM \dot{ \bfr}^n$};
			\draw[-stealth] (2.95,4.3) -- node[right, pos=0.35] {\tiny $\quad \partial^\tau_q$} (3.45,4.6+\yshftii);
			
			\pgfmathsetmacro{\xshfti}{-0.83}
			\draw[-stealth] (1.4+\xshfti,5.05) -- (-0.675+\xshfti,5.05) -- node[below, midway, sloped] {\footnotesize \hspace*{5 mm} test with $\bfe_u^n-\eta \bfe_u^{n-1}$} (-0.675+\xshfti,0.5);
			\draw[-stealth] (1.4+\xshfti,4.55) -- (-0.525+\xshfti,4.55) -- node[above, midway, sloped] {\footnotesize \hspace*{2.5 mm} test with $\bfe_w^n$} (-0.525+\xshfti,0.5);
			\draw (-0.6+\xshfti,0.4) node {\tiny $+$};
			
			\draw[-stealth] (1.4+\xshfti,4.95) --  (-2.075+\xshfti,4.95) -- node[below, midway, sloped] {\footnotesize \hspace*{7.5 mm} test with $\bfe_w^n-\eta \bfe_w^{n-1}$} (-2.075+\xshfti,0.5);
			\draw[-stealth] (1.4+\xshfti,4.45) -- node[below, pos=1.1+\xshfti, sloped] {\tiny \hspace*{1 mm} $(\cdot)^n-\eta (\cdot)^{n-1}$} (-1.925+\xshfti,4.45) -- node[above, midway, sloped] {\footnotesize \hspace*{2.5 mm} test with $\dot{\bfe}_u^n$}  (-1.925+\xshfti,0.5);
			\draw (-2+\xshfti,0.4) node {\tiny $-$};
			
			\draw (-1.3+\xshfti,5.35) node {\tiny (i)};
			
			\pgfmathsetmacro{\xshftii}{2.18}
			\pgfmathsetmacro{\xshftiib}{-0.6}
			\draw[-stealth] (2.2+\xshftii,5.05) -- (7.925,5.05) -- node[below, midway, sloped] {\footnotesize \hspace*{15 mm} test with $\dot{\bfe}_u^n$} (7.925,0.5);
			\draw[-stealth] (5.3,4.55+\yshftii) -- (8.075,4.55+\yshftii) -- node[above, midway, sloped] {\footnotesize  test with $\bfe_w^n-\eta \bfe_w^{n-1}$} (8.075,0.5);
			\draw (8,0.4) node {\tiny $+$};
			
			\draw[-stealth] (2.2+\xshftii,4.95) -- node[below, midway, sloped] {\tiny  $(\cdot)^n-\eta (\cdot)^{n-1}$}(9.925+\xshftiib,4.95) -- node[below, midway, sloped] {\footnotesize \hspace*{15.5 mm} test with $\dot{\bfe}_w^n$} (9.925+\xshftiib,0.5);
			\draw[-stealth] (5.3,4.45+\yshftii) -- (10.075+\xshftiib,4.45+\yshftii) -- node[above, midway, sloped] {\footnotesize \hspace*{0 mm} test with $\dot{\bfe}_u^n -\eta \dot{\bfe}_u^{n-1}$} (10.075+\xshftiib,0.5);
			\draw (10+\xshftiib,0.4) node {\tiny $-$};
			
			\draw (9.0+\xshftiib/2,5.35) node {\tiny (ii)};

			\draw[draw, align=left] (4.15,-1.1) node {Stability estimate: \\ Proposition~\ref{Prop:StabilityEst}};
			\draw[draw, align=left] (2.5,-0.1) -- node[above, midway, sloped] {\footnotesize combining (i) and (ii)}  (6,-0.1);
			\draw[-stealth] (4.15,-0.1) -- (4.15,-0.65);
		\end{tikzpicture}
	\end{adjustbox}
	\caption{Sketch of the energy estimates for the fully discrete stability proof.} \label{fig:Scheme}
\end{figure}

In the following $c$ and $C$ are generic constants, independent of $h$ and $\tau$, that may take different values on different occurrences. By $\rho, \mu, \widetilde \mu >0$ we denote a small number, independent of $h$ and $\tau$, used in Young's inequalities, and hence we often incorporate multiplicative constants into those yet unchosen factors.

\textit{(Preparations)} 
In Section~\ref{section:Consistency} it will be shown that the defects are bounded by $C(h^2 + \tau^q)$, hence, under the mild step-size restriction $\tau^q \leq C_0 h^\kappa$, the assumed bounds in \eqref{eq:assumed bounds defects} are all satisfied.

By our assumptions $\bfW'$ and $\bfW''$ are \emph{locally} Lipschitz continuous. Since the extrapolations in the non-linear terms only depend on past values, we use a bootstrapping argument to show that $\mathbf{\widetilde u}^n$ is $L^\infty$-close to $\mathbf{\widetilde u}_\ast^n$, as indicated in Section~\ref{subsec:Linftyu}.

Let $ t_{\max} \in (0,T]$ be the maximal time such that the following inequality holds for $k \tau \leq (n-1)\tau \leq t_{\max}$:
\begin{equation}
\label{eq:Linftybounde_u}
	\| \gamma_h \bfe^{k}_u \|_{L^\infty(\Ga_h)} \leq \| \bfe^{k}_u \|_{L^\infty(\Om_h)} 
	\leq h^{\frac{\kappa}{2}-\frac{d}{4}}. 
\end{equation}
Note that the first inequality holds in general for bulk--surface finite element functions (since $\partial \Om_h = \Ga_h$), and for the starting values this bound is directly satisfied by \eqref{eq:defI}, for $h$ sufficiently small. From now on we assume in the proof that $(n-1)\tau \leq t_{\max}$ and therefore these bounds hold. At the end of the proof we will then show that $t_{\max}$ coincides with $T$. 

By \eqref{eq:Linftybounde_u} it immediately follows that, for any $k = q+1,\dotsc,n$, 
\begin{align*}
	\| \gamma_h (\mathbf{\widetilde u}^k-\mathbf{\widetilde u}_\ast^k) \|_{L^\infty(\Ga_h)} \leq&\ \| \mathbf{\widetilde u}^k-\mathbf{\widetilde u}_\ast^k \|_{L^\infty(\Om_h)} \leq \sum_{j=0}^{q-1} |\gamma_j|  \|\bfe_u^{k-1-j}\|_{L^\infty} \leq c h^{\frac{\kappa}{2}-\frac{d}{4}}.
\end{align*}
Together with \eqref{eq:exactuLinfty}, we can conclude that there is a large enough compact set, which contains both $\mathbf{\widetilde u}^k$ and $\mathbf{\widetilde u}^k_\ast$. On this set $\bfW'$ and $\bfW''$ are Lipschitz continuous.

\textit{(i)} 
We will now perform the energy estimates sketched in the branch on the left-hand side of Figure~\ref{fig:Scheme}. Let $q+1 \leq k \leq n$, such that $(n-1) \tau \leq t_{\max}$. Following the scheme in Figure~\ref{fig:Scheme}, we start by testing \eqref{eq:CHNumericalScheme1} with $\bfe^k_u - \eta \bfe^{k-1}_u$, where the multiplier factor $\eta$ is chosen for a given $q$ as in Lemma~\ref{Lemma:MultiplierTechnique}, and adding \eqref{eq:CHNumericalScheme2}, tested with $\bfe^k_w$, to obtain
\begin{align}
	(\bfde^k_u) \TbfM (\bfe^k_u - \eta \bfe^{k-1}_u) + (\bfe^k_w) \TbfM \bfe^k_w - \eta (\bfe^k_w) \TbfA \bfe^{k-1}_u = (\bfe^k_w) \TbfM \bfr^k - (\bfe^k_u-\eta \bfe^{k-1}_u) \TbfM \bfd^k_u. \label{eq:i1}
\end{align}

Next we consider \eqref{eq:CHNumericalScheme2} for $k$, subtract $\eta$ times \eqref{eq:CHNumericalScheme2} for $k-1$, and test this equation with $\bfde^k_u$. Then we subtract \eqref{eq:CHNumericalScheme1}, tested with $(\bfe^k_w - \eta \bfe^{k-1}_w)$, to obtain
\begin{align}
	(\bfde^k_u) \TbfA (\bfe^k_u - \eta \bfe^{k-1}_u) + (\bfe^k_w) \TbfA \bfe^k_w - \eta (\bfe^k_w) \TbfA \bfe^{k-1}_w =&\   - (\bfe^k_w - \eta \bfe^{k-1}_w) \TbfM \bfd^k_u - (\bfde^k_u) \TbfM (\bfr^k - \eta \bfr^{k-1}). \label{eq:i2}
\end{align}

Adding up the equalities \eqref{eq:i1} and \eqref{eq:i2}, recalling that $\bfK = \bfM + \bfA$, we obtain
\begin{align}
	(\bfde^k_u) \TbfK (\bfe^k_u - \eta \bfe^{k-1}_u) + \| \bfe^k_w\|_\bfK^2 =&\  \eta (\bfe^k_w) \TbfA \bfe^{k-1}_w + \eta (\bfe^k_w) \TbfA \bfe^{k-1}_u + (\bfe^k_w) \TbfM \bfr^k \nonumber \\
	&\ - (\bfe^k_u-\eta \bfe^{k-1}_u) \TbfM \bfd^k_u - (\bfe^k_w - \eta \bfe^{k-1}_w) \TbfM \bfd^k_u \nonumber \\
	&\ - 	(\bfde^k_u) \TbfM (\bfr^k - \eta \bfr^{k-1}). \label{eq:equalityi}
\end{align}

We start by estimating the first term on the left-hand side. We apply Lemma~\ref{Lemma:DahlquistOrig} and Lemma~\ref{Lemma:MultiplierTechnique} to obtain, with the symmetric positive definite matrix $G$ and $\bfE^k_u=(\bfe^k_u,\dotsc,\bfe^{k-q+1}_u)$:
\begin{align}
	\frac{1}{\tau} \big(\|\bfE^k_u\|_{G,\bfK}^2 - \|\bfE^{k-1}_u\|_{G,\bfK}^2\big) \leq (\bfde^k_u) \TbfK  (\bfe^k_u - \eta \bfe^{k-1}_u). \label{eq:idahlquist}
\end{align}

Therefore, we have, with \eqref{eq:equalityi} and \eqref{eq:idahlquist}, multiplying by $\tau$, and summing over $k=q+1,\dotsc,n$, that
\begin{align}
	\big(\|\bfE^n_u\|_{G,\bfK}^2 - \|\bfE^{q}_u\|_{G,\bfK}^2\big) + \tau \sum_{k=q+1}^n \| \bfe^k_w\|_\bfK^2 \leq&\ \eta \tau \sum_{k=q+1}^n \bigg[ (\bfe^k_w) \TbfA \bfe^{k-1}_w + (\bfe^k_w) \TbfA \bfe^{k-1}_u \bigg] \nonumber \\
	&\ - \tau \sum_{k=q+1}^n \bigg[ (\bfe^k_u-\eta \bfe^{k-1}_u) \TbfM \bfd^k_u + (\bfe^k_w - \eta \bfe^{k-1}_w) \TbfM \bfd^k_u \bigg] \nonumber \\
	&\ + \tau \sum_{k=q+1}^n \bigg[ (\bfe^k_w) \TbfM \bfr^k - (\bfde^k_u) \TbfM (\bfr^k - \eta \bfr^{k-1}) \bigg] \nonumber \\
	=:&\ I + II + III. \label{eq:iLHS}
\end{align}

The terms on the right-hand side of \eqref{eq:iLHS} are estimated separately, using the Cauchy--Schwarz inequality and Young's inequality. For the first term we have
\begin{align}
	I \leq&\ \frac{\eta}{2} \tau \sum_{k=q+1}^n \bigg(\|\bfe_w^k\|_\bfA^2 + \|\bfe_w^{k-1}\|_\bfA^2 \bigg) + \rho \tau \sum_{k=q+1}^n \|\bfe_w^k\|_\bfA^2 + c \tau \sum_{k=q+1}^n \|\bfe_u^{k-1}\|_\bfA^2 \nonumber \\
	\leq&\ (\eta + \rho) \tau \sum_{k=q+1}^n \| \bfe_w^k \|_\bfK^2 + c\tau \sum_{k=q+1}^{n-1} \|\bfe_u^k\|_\bfK^2 + \frac{\eta}{2} \tau \| \bfe_w^q \|_\bfK^2 + c \tau \|\bfe_u^q \|_\bfK^2 .\label{eq:iRHS1}
\end{align}

For the second term on the right-hand side of \eqref{eq:iLHS} we use the definition of the dual norm and Young's inequality to obtain
\begin{align}
	II \leq&\ \rho \tau \sum_{k=q+1}^n \bigg(\| \bfe_w^k\|_\bfK^2 + \| \bfe_w^{k-1}\|_\bfK^2 \bigg) + \mu \tau \sum_{k=q+1}^n \bigg(\| \bfe_u^k\|_\bfK^2 + \| \bfe_u^{k-1}\|_\bfK^2 \bigg) + c \tau \sum_{k=q+1}^n \| \bfd_u^k\|_{\ast,h}^2 \nonumber \\
	\leq&\ \rho \tau \sum_{k=q+1}^n \| \bfe_w^k\|_\bfK^2 + \mu \tau \sum_{k=q+1}^n \| \bfe_u^k\|_\bfK^2 + \rho \tau \| \bfe_w^q \|_\bfK^2 + \mu \tau \| \bfe_u^q \|_\bfK^2 + c \tau \sum_{k=q+1}^n \| \bfd_u^k\|_{\ast,h}^2 . \label{eq:iRHS2}
\end{align}

For the third term $III$ on the right-hand side of \eqref{eq:iLHS} we recall that $\bfr^k= \bfW'(\mathbf{\widetilde u}^k) - \bfW'(\mathbf{\widetilde u}^k_\ast) - \bfd^k_w + \bftheta^k$, see~\eqref{eq:defr}. The defect terms are estimated by Young's inequality and the definition of the shift \eqref{eq:defShift}:
\begin{align}
	\tau \sum_{k=q+1}^n &\bigg[ -(\bfe^k_w) \TbfM (\bfd_w^k-\bftheta^k) + (\bfde^k_u) \TbfM \Big(\bfd_w^k - \eta \bfd_w^{k-1}-(\bftheta^k-\eta\bftheta^{k-1})\Big) \bigg] \nonumber \\
	\leq&\ \rho \tau \sum_{k=q+1}^n \| \bfe_w^k \|_\bfK^2 + \widetilde \mu \tau \sum_{k=q+1}^n \|\bfde_u^k\|_\bfK^2 + c \tau \sum_{k=q}^n \|\bfd_w^k \|_{\ast,h}^2 + c \| \bfd_w^{q-1} \|_{\ast,h}^2. \label{eq:iRHS3}
\end{align}

On the other hand, for the non-linear terms we use the Lipschitz continuity of $W'$ in a compact set, as explained in the preparations above. 
%
%
%
%
%
We estimate by 
\begin{align}
	\tau \sum_{k=q+1}^n &\bigg[ (\bfe^k_w - \bfde^k_u) \TbfM \big(\bfW'(\widetilde \bfu^k) - \bfW'(\widetilde \bfu_\ast^k)\big) + \eta (\bfde^k_u) \TbfM \big(\bfW'(\widetilde \bfu^{k-1}) - \bfW'(\widetilde \bfu_\ast^{k-1})\big)  \bigg] \nonumber\\
	\leq&\ \rho \tau \sum_{k=q+1}^n \| \bfe_w^k\|_\bfK^2 + \widetilde \mu \tau \sum_{k=q+1}^n \| \bfde_u^k\|_\bfK^2 + c \tau \sum_{k=q}^n \| \widetilde \bfe_u^k \|_\bfK^2 \nonumber\\
	\leq&\ \rho \tau \sum_{k=q+1}^n \| \bfe_w^k\|_\bfK^2 + \widetilde \mu \tau \sum_{k=q+1}^n \| \bfde_u^k\|_\bfK^2 + c \tau \sum_{k=q+1}^{n-1} \| \bfe_u^k \|_\bfK^2 + c \tau \| \bfe_u^q \|_\bfK^2 + c I_h^{q-1}, \label{eq:iRHS4}
\end{align}
where $I_h^{q-1}$ is the error in the initial values, see \eqref{eq:defI}.

Combining inequality \eqref{eq:iLHS} with the bounds \eqref{eq:iRHS1}, \eqref{eq:iRHS2}, \eqref{eq:iRHS3}, and \eqref{eq:iRHS4}, 
using an absorption by choosing $\rho >0$ such that $0<\eta + \rho <1$, 
%
applying the estimates on the $G$-norm by the $\|\cdot\|_{\bfK}$-norm \eqref{eq:GnormEst}, and after a final absorption of the terms involving $\mu$, we obtain
\begin{align*}
	\|\bfe^n_u\|_\bfK^2 + \tau \sum_{k=q+1}^n \| \bfe^k_w\|_\bfK^2 \leq&\  \widetilde \mu \tau \sum_{k=q+1}^n \|\bfde_u^k\|_\bfK^2 + c\tau \sum_{k=q+1}^{n-1} \|\bfe_u^k\|_\bfK^2 + c \tau \| \bfe_w^q \|_\bfK^2 + c \tau \|\bfe_u^q \|_\bfK^2 \\
	&\ + c \tau \sum_{k=q+1}^n \| \bfd_u^k\|_{\ast,h}^2 + c \tau \sum_{k=q}^n \|\bfd_w^k \|_{\ast,h}^2 + c \| \bfd_w^{q-1} \|_{\ast,h}^2 + cI_h^{q-1} .
\end{align*}
The only terms which we need to work on further are the terms involving $\bfe^q_{u}$ and $\bfe^q_w$, here we insert the estimate \eqref{eq:qest} for the special case $n=q$ and arrive at the first energy estimate 
\begin{equation}
	\label{eq:energyi}
	\begin{aligned}
		\|\bfe^n_u\|_\bfK^2 + \tau \sum_{k=q+1}^n \| \bfe^k_w\|_\bfK^2 
		\leq &\ \widetilde \mu \tau \sum_{k=q+1}^n \| \bfde^k_u \|_\bfK^2 + c \tau \sum_{k=q+1}^{n-1} \|\bfe^k_u \|_\bfK^2 + c \big( I_h^{q-1} + D_h^n \big) ,
	\end{aligned}
\end{equation}
where we recall the notations \eqref{eq:defI} and \eqref{eq:defD}.

\textit{(The special case $n=q$)} For $n=q$ it remains to show an estimate analogous to \eqref{eq:StabilityEst}. By slightly modifying the argument for the general case we obtain
\begin{align}
\label{eq:qest}
	\| \bfe^q_u \|_\bfK^2 + \| \bfe^q_w \|_\bfK^2 + \tau \| \bfde^q_u \|_\bfK^2  \leq&\  C I_h^{q-1} + C D_h^q .
\end{align}	
The difference to the general case is that there is no equation for $n-1=q-1$. We solve this problem by adding the missing terms on both sides, so that we can use Dahlquist's $G$-stability result (Lemma~\ref{Lemma:DahlquistOrig} and \ref{Lemma:Dahlquist}) and the multiplier technique of Nevanlinna and Odeh (Lemma~\ref{Lemma:MultiplierTechnique}), these terms involve only initial values and therefore we estimate them directly.

%

\textit{(ii)} 
We will now perform the energy estimates sketched in the branch on the right-hand side of Figure~\ref{fig:Scheme}. Let again $q+1 \leq k \leq n$, such that $(n-1) \tau \leq t_{\max}$. For the second energy estimate, we take the discrete time derivative of the second equation. To this end, for $k\geq 2q$ we consider the linear combination of \eqref{eq:CHNumericalScheme2}, for $k-i=k-q,\dotsc,k$, weighted by $\delta_i / \tau$, which gives
\begin{align}
	\bfM \bfde^k_w - \bfA \bfde^k_u = \bfM \mathbf{\dot r}^k , \qquad k\geq 2q. \label{eq:difCH2}
\end{align}

For $k < 2q$ we have to use a different equation, since then in the discrete derivative the initial values are starting to appear, and for the initial values we do not have the equation \eqref{eq:CHNumericalScheme2}. So instead we have to add the terms involving the initial values on both sides. For convenience of notation we define 
\begin{align*}
	\bfM \bfr^i:= \bfM \bfe^i_w - \bfA \bfe^i_u = -\bfM \bfd_w^i + \bfM \bftheta^i, \qquad \text{for} \quad i = 0, \dotsc, q-1 ,
\end{align*}
where the second equality exactly holds by the definition of the initial values and our definition of the shift $\bftheta^n$, see Section~\ref{subsec:modifiedEq} and especially \eqref{eq:defShift}. This definition allows us to write the equation identical to \eqref{eq:difCH2} even for $q+1\leq k <2q$. We now test \eqref{eq:CHNumericalScheme1} and \eqref{eq:difCH2} in a similar way as before.

Testing \eqref{eq:CHNumericalScheme1} with $\bfde^k_u$ and adding \eqref{eq:difCH2}, tested with $\bfe^k_w - \eta \bfe^{n-1}_w$, yields
\begin{align}
	(\bfde^k_u) \TbfM \bfde^k_u + (\bfe^k_w - \eta \bfe^{k-1}_w) \TbfM \bfde^k_w + \eta (\bfe^{k-1}_w) \TbfA \bfde^k_u = (\bfe^k_w - \eta \bfe^{k-1}_w) \TbfM \mathbf{\dot r}^k - (\bfde^k_u) \TbfM \bfd^k_u. \label{eq:ii.1}
\end{align}
Using again the theory of Dahlquist and Nevanlinna \& Odeh, i.e., Lemma~\ref{Lemma:DahlquistOrig} and Lemma~\ref{Lemma:MultiplierTechnique}, we estimate
\begin{align*}
	| \bfde^k_u |_\bfM^2 + \frac{1}{\tau} \Big( \|\bfE^k_w \|_{G,\bfM}^2 - \| \bfE^{k-1}_w\|_{G,\bfM}^2 \Big) \leq (\bfe^k_w - \eta \bfe^{k-1}_w) \TbfM \mathbf{\dot r}^k - (\bfde^k_u) \TbfM \bfd^k_u - \eta (\bfe^{k-1}_w) \TbfA \bfde^k_u, 
\end{align*} 
and sum over $k=q+1,\dotsc,n$, multiply by $\tau$, and rearrange the terms to obtain 
\begin{align}
	\Big( \|\bfE^n_w\|_{G,\bfM}^2 - \|\bfE^q_w\|_{G,\bfM}^2 \Big)  + \tau \sum_{k=q+1}^n \|\bfde^k_u \|_\bfM^2  
	\leq &\   - \tau \sum_{k=q+1}^n (\bfde^k_u) \TbfM \bfd^k_u - \eta \tau \sum_{k=q+1}^n (\bfe^{k-1}_w) \TbfA \bfde^k_u \nonumber \\
	&\ + \tau \sum_{k=q+1}^n (\bfe^k_w - \eta \bfe^{k-1}_w) \TbfM \mathbf{\dot r}^k \nonumber \\
	=:&\ I + II + III. \label{eq:ii.1LHS}
\end{align}
As in Part (i) of the proof, we use Cauchy--Schwarz inequality and Young's inequality to estimate the three terms separately. For the first two terms, similarly as before, we obtain
\begin{align}
	I \leq &\ \widetilde \mu \tau \sum_{k=q+1}^n \|\bfde^k_u\|_\bfK^2 + c \tau \sum_{k=q+1}^n \|\bfd^k_u \|_{\ast,h}^2, \label{eq:ii.1rhs1} \\
	II	\leq &\ \widetilde \mu \tau \sum_{k=q+1}^n \|\bfde^k_u\|_\bfK^2 + c \tau \sum_{k=q}^{n-1} \|\bfe^k_w \|_\bfK^2.\label{eq:ii.1rhs2}
\end{align}

In order to analyse the last term $III$ on the right-hand side of \eqref{eq:ii.1LHS}, we analyse the term tested by $\bfe^k_w$. The term tested by $\bfe^{k-1}_w $ can be treated analogously. 

We start by splitting $\bfr^j$ into two parts, and define	
\begin{align*}
	\bfM \bfr_1^j := &\ \begin{cases}
		\bfM \bfW'(\widetilde \bfu^j) - \bfM \bfW'(\widetilde \bfu_\ast^j) , \quad & \text{for } \ j\geq q , \\
		0 , \quad & \text{for } \  0 \leq j < q ;
	\end{cases} ,\\
	\bfM \bfr_2^j:=&\ -\bfM \bfd_w^j + \bfM \bftheta^j .
\end{align*}

In order to estimate the terms involving $\bfr_1$, according to the proof of Lemma~\ref{Lemma:partialtodot}, we factorize the generating polynomial into $\delta(\zeta) = (1-\zeta) \sigma(\zeta)$, and, upon recalling \eqref{eq:factorization}, we obtain 
\begin{align*}
	\tau \sum_{k=q+1}^n (\bfe^k_w) \TbfM \mathbf{\dot r}^k =&\ \tau \sum_{k=q+1}^n \sum_{j=k-q+1}^k \sigma_{k-j} (\bfe^k_w) \TbfM \partial^\tau \bfr^j .
\end{align*}
This allows us to write the discrete derivative as a linear combination of difference quotients. We use the Cauchy--Schwarz inequality, \emph{local Lipschitz continuity} of $W''$ (cf.~preparations), and the definition $\mathbf{\widetilde u}^j_\theta = \mathbf{\widetilde u}^j_\ast + \theta (\mathbf{\widetilde u}^j-\mathbf{\widetilde u}^j_\ast)$, for $0\leq \theta \leq 1$, to estimate for $q+1\leq j \leq n$: 
\begin{align}
	(\bfe^k_w) \TbfM \partial^\tau \big(\bfW'(\mathbf{\widetilde u}^j)-\bfW'(\mathbf{\widetilde u}^j_\ast)\big) =&\ (\bfe^k_w) \TbfM \int_0^1 \frac{\d}{\d\theta} \partial^\tau \bfW'(\mathbf{\widetilde u}^j_\theta) ~\d \theta \nonumber \\
	=&\ (\bfe^k_w) \TbfM \Big(\int_0^1 \partial^\tau \big(\bfW''(\mathbf{\widetilde u}^j_\theta) \mathbf{\widetilde e}^j_u\big) ~\d \theta \Big)\nonumber \\
	=&\ (\bfe^k_w) \TbfM \Big(\int_0^1 \bfW''(\mathbf{\widetilde u}^j_\theta) \partial^\tau \mathbf{\widetilde e}^j_u ~\d \theta + \int_0^1 \big(\partial^\tau \mathbf{W''}(\mathbf{\widetilde u}^j_\theta)\big) \mathbf{\widetilde e}^{j-1}_u ~\d \theta\Big) \nonumber \\
	\leq&\  c \, \|\bfe^k_w\|_\bfM \Big(\|\partial^\tau \mathbf{\widetilde e}^j_u\|_\bfM + \|\partial^\tau \mathbf{\widetilde u}^j_\ast\|_{L^\infty(\Om_h)} \|\mathbf{\widetilde e}^{j-1}_u\|_\bfM + \| \partial^\tau \mathbf{\widetilde e}^j_u \|_\bfM \| \mathbf{\widetilde e}^{j-1}_u\|_{L^\infty(\Om_h)} \Big) \nonumber \\
	\leq&\  c \, \|\bfe^k_w\|_\bfM \bigg( \sum_{l=j-q}^{j-1} \|\partial^\tau \bfe^l_u \|_\bfM + \sum_{l=j-q}^{j-1} \| \bfe^l_u \|_\bfM\bigg), 
	\label{eq:r1.1}
\end{align}
%
where, in the last inequality, we used the following $L^\infty$-bounds, cf.~\eqref{eq:exactuLinfty} and \eqref{eq:Linftybounde_u}:
\begin{align*}
	\| \mathbf{\widetilde e}^{j-1}_u\|_{L^\infty(\Om_h)} \leq&\  c \sum_{i=j-1-q}^{j-1} \| \bfe^{i}_u\|_{L^\infty(\Om_h)} \leq c, \\
	\|\partial^\tau \mathbf{\widetilde u}^j_\ast\|_{L^\infty(\Om_h)} 	\leq&\ c \| \bfu_\ast^{(1)} \|_{L^\infty([0,T],L^\infty(\Om_h))} \leq c.
\end{align*} 

By our choice of initial value shifts, for $j\leq q$ we have that $\mathbf{\widetilde u}^q = \mathbf{\widetilde u}_\ast^q$, and hence $\bfr_1^q=0$ and so $\partial^\tau \bfr_1^j=0$. 
%
Therefore, by \eqref{eq:r1.1}, Lemma~\ref{Lemma:partialtodot} (which estimates $\partial^\tau$ in terms of $\partial^\tau_q$), and the estimate \eqref{eq:qest} for the special case $n=q$, we have that
\begin{align}
	\tau \sum_{k=q+1}^n \sum_{j=k-q+1}^k \sigma_{k-j} (\bfe^k_w) \TbfM \partial^\tau \bfr_1^j \leq&\  \widetilde \mu \tau \sum_{k=q+1}^{n-1} \| \bfde_u^k \|_\bfM^2 +  \mu \tau \sum_{k=q+1}^{n-1} \| \bfe_u^k \|_\bfM^2 + c \tau \sum_{k=q+1}^n \| \bfe_w^k \|_\bfM^2 +cI_h^{q-1} + cD_h^q. \label{eq:non-linear1}
\end{align}

In order to estimate the terms involving $\bfr_2$, we recall the definition of the shift \eqref{eq:defShift}, and therefore have
\begin{align}
	\partial^\tau \bftheta^j= \begin{cases}
		0; \qquad &q\leq j\\
		\partial^\tau \bfd_w^j; \qquad &1\leq j \leq q-1
	\end{cases}. \label{eq:r2.2}
\end{align}
Using now \eqref{eq:r2.2} and once again the factorization \eqref{eq:factorization} of the generating polynomial $\delta(\zeta) = (1-\zeta) \sigma(\zeta)$, we have 
\begin{align}
	\boldsymbol{\dot \vartheta}^k= \sum_{i=k-q+1}^{q-1} \sigma_i \partial^\tau \bfd_w^i. \label{eq:dShift}
\end{align} 
Then, we estimate using the definition of the dual norm and Young's inequality:
\begin{align}
	\tau \sum_{k=q+1}^n (\bfe_w^k) \TbfM \mathbf{\dot r}_2^k =&\ \tau \sum_{k=q+1}^n (\bfe_w^k) \TbfM \Big(\mathbf{\dot d}_w^k + \boldsymbol{\dot \vartheta}^k\Big) \nonumber \\
	\leq&\ \rho \tau \sum_{k=q+1}^n \|\bfe_w^k\|_\bfK^2 + c \tau \sum_{k=q+1}^n \|\mathbf{\dot d}_w^k\|_{\ast,h}^2 + c \tau \sum_{k=q+1}^n \|\boldsymbol{\dot \vartheta}^k\|_{\ast,h}^2 \nonumber \\
	\leq&\ \rho \tau \sum_{k=q+1}^n \|\bfe_w^k\|_\bfK^2 + c \tau \sum_{k=q+1}^n \|\mathbf{\dot d}_w^k\|_{\ast,h}^2 + c \tau \sum_{i=1}^{q-1} \| \partial^\tau \bfd_w^i\|_{\ast,h}^2\label{eq:non-linear2}
\end{align}

In order to obtain the bound on the full term, we add \eqref{eq:non-linear1} and \eqref{eq:non-linear2} and use an analogous estimate for $\eta \bfe_w^{k-1}$, to arrive at
\begin{align}
	III=  \tau \sum_{k=q+1}^n (\bfe_w^k- \eta \bfe_w^{k-1})\TbfM \mathbf{\dot r}^k \leq&\ \widetilde \mu \tau \sum_{k=q+1}^{n-1} \| \bfde_u^k \|_\bfM^2 + \mu \tau \sum_{k=q+1}^{n-1} \| \bfe_u^k \|_\bfM^2 + c \tau \sum_{k=q+1}^n \| \bfe_w^k \|_\bfM^2  \nonumber \\
	&\ +\rho \tau \sum_{k=q+1}^n \|\bfe_w^k\|_\bfK^2 
	+ c\big(I_h^{q-1} + D_h^n\big). \label{eq:rdotest}
\end{align}

Then we plug \eqref{eq:ii.1rhs1}, \eqref{eq:ii.1rhs2}, and \eqref{eq:rdotest} into \eqref{eq:ii.1LHS}, and use the estimates from \eqref{eq:GnormEst} for the $G$-norm, as well as the estimate \eqref{eq:qest}, to obtain the first part of the second energy estimate:
\begin{align}
	\tau \sum_{k=q+1}^n \|\bfde^k_u \|_\bfM^2 + \sum_{k=n-q+1}^n\|\bfe^k_w\|_\bfM^2 \leq&\  \widetilde \mu \tau \sum_{k=q+1}^n \| \bfde^k_u \|_\bfK^2  + c \tau \sum_{k=q+1}^{n} \| \bfe^k_w \|_\bfK^2 + c \tau \sum_{k=q+1}^{n-1}  \| \bfe^k_u\|_\bfK^2  + c \big( I_h^{q-1} + D_h^n  \big) . \label{eq:ii.1bound}
\end{align}

For the second branch of the right-hand side of Figure~\ref{fig:Scheme}, for $q+1 \leq k \leq n$, we consider \eqref{eq:CHNumericalScheme1} and subtract $\eta$ times the same equation for $k-1$,   and test with $\bfde^k_w$, then subtract \eqref{eq:difCH2} tested with $\bfde^k_u - \eta \bfde^{k-1}_u$, to obtain
\begin{align*}
	(\bfde^k_w) \TbfA (\bfe^k_w - \eta \bfe^{k-1}_w) + (\bfde^k_u) \TbfA \bfde^k_u -\eta (\bfde^{k-1}_u) \TbfA \bfde^k_u =&\ - (\bfde^k_w) \TbfM (\bfd^k_u - \eta \bfd^{k-1}_u) - (\bfde^k_u - \eta \bfde^{k-1}_u) \TbfM \mathbf{\dot r}^k.
\end{align*}

Using again the theory of Dahlquist and Nevanlinna \& Odeh, Lemma~\ref{Lemma:Dahlquist} and Lemma~\ref{Lemma:MultiplierTechnique}, we estimate, as before, by
\begin{align*}
	\frac{1}{\tau} \Big( \|\bfE^k_w\|_{G,\bfA}^2 - \|\bfE^{k-1}_w\|_{G,\bfA}^2 \Big) + \| \bfde^k_u \|^2_\bfA \leq&\  \eta (\bfde^{k-1}_u) \TbfA \bfde^k_u - (\bfde^k_w) \TbfM (\bfd^k_u - \eta \bfd^{k-1}_u) - (\bfde^k_u - \eta \bfde^{k-1}_u) \TbfM \mathbf{\dot r}^k.
\end{align*}
Now we sum this inequality from $k=q+1$ to $n$ and multiply by $\tau$, which yields
\begin{align}
	\Big( \|\bfE^n_w\|_{G,\bfA}^2 - \|\bfE^q_w\|_{G,\bfA}^2 \Big) + \tau \sum_{k=q+1}^n \|\bfde^k_u \|^2_\bfA \leq&\  \tau \sum_{k=q+1}^n \eta (\bfde^{k-1}_u) \TbfA \bfde^k_u - \tau \sum_{k=q+1}^n (\bfde^k_u - \eta \bfde^{k-1}_u) \TbfM \mathbf{\dot r}^k \nonumber \\
	&\ - \tau \sum_{k=q+1}^n (\bfde^k_w) \TbfM (\bfd^k_u - \eta \bfd^{k-1}_u) \nonumber \\
	=:&\ I + II + III . \label{eq:ii.2.lhs}
\end{align}	

The first term on the right-hand side is estimated by Cauchy-Schwarz and Young's inequality, and we obtain  
\begin{align}
	I \leq&\  \frac{\eta}{2} \tau \sum_{k=q}^{n-1} \| \bfde^{k}_u \|^2_\bfA +  \frac{\eta}{2} \tau \sum_{k=q+1}^{n} \| \bfde^{k}_u \|^2_\bfA \leq  \eta \tau \sum_{k=q}^{n} \| \bfde^{k}_u \|^2_\bfA. \label{eq:ii.2.rhs2}
\end{align}

The second term on the right-hand side we estimate as in \eqref{eq:rdotest}, just with $\bfde^k_u$ instead of $\bfe^k_w$, and we obtain
\begin{align}
	II \leq&\  c \tau \sum_{k=q+1}^n \|\bfde^k_u \|_\bfM^2 + \widetilde \mu \tau \sum_{k=q+1}^n \| \bfde^k_u \|_\bfK^2 + \mu \tau \sum_{k=q+1}^{n-1} \| \bfe^k_u \|_\bfK^2 
	+c\big(I_h^{q-1} + D_h^n\big). \label{eq:ii.2.rhs1}
\end{align}


In order to estimate the last term on the right-hand side we have to overcome the problem that we do not have a good control on $\bfde^k_w$. We therefore use a discrete product rule to transfer the discrete derivative onto $\bfd^k_u$, cf.~\cite[Proposition~10.1]{KovacsLiLubich2019} \cite[Proposition~5.1]{HarderKovacs2022}. We aim at showing
\begin{align}
	III= \tau \sum_{k=q+1}^n (\bfde^k_w) \TbfM (\bfd^k_u - \eta \bfd^{k-1}_u) \leq  \rho \sum_{\substack{k=n-q+1 \\ k\geq q+1}}^{n} \|\bfe^{k}_w \|_\bfK^2 + \rho \tau \sum_{k=q+1}^n \| \bfe^k_w \|_\bfK^2+c\big(I_h^{q-1} + D_h^n\big) . \label{eq:boundproductrule}
\end{align}
Using the factorization $\delta(\zeta) = (1-\zeta) \sigma(\zeta)$ of Lemma~\ref{Lemma:partialtodot} once again, we obtain
\begin{align*}
	\tau \sum_{k=q+1}^n (\bfde^k_w) \TbfM (\bfd^k_u- \eta \bfd^{k-1}_u) = \tau \sum_{k=q+1}^n \sum_{j=0}^{q-1} \sigma_j (\partial^\tau \bfe^{k-j}_w) \TbfM (\bfd^k_u- \eta \bfd^{k-1}_u).
\end{align*} 
Then, using the identities
\begin{align*}
	(\partial^\tau \bfe^{k-j}_w) \TbfM \bfd^k_u =&\ \partial^\tau \big((\bfe^{k-j}_w) \TbfM \bfd^k_u\big) - (\bfe^{k-j-1}_w) \TbfM \partial^\tau \bfd^k_u, \\
	(\partial^\tau \bfe^{k-j}_w) \TbfM \bfd^{k-1}_u =&\ \partial^\tau \big((\bfe^{k-j}_w) \TbfM \bfd^k_u\big) - (\bfe^{k-j}_w) \TbfM \partial^\tau \bfd^k_u,
\end{align*}
we arrive at
\begin{align*}
	\tau \sum_{k=q+1}^n (\bfde^k_w) \TbfM (\bfd^k_u - \eta \bfd^{k-1}_u) =&\ \tau (1-\eta)\sum_{k=q+1}^n \sum_{j=0}^{q-1} \sigma_j \partial^\tau \big((\bfe^{k-j}_w) \TbfM \bfd^k_u\big) - \tau \sum_{k=q+1}^n \sum_{j=0}^{q-1} \sigma_j (\bfe^{k-j-1}_w) \TbfM \partial^\tau \bfd^k_u \nonumber \\
	&\ + \tau \eta \sum_{k=q+1}^n \sum_{j=0}^{q-1} \sigma_j (\bfe^{k-j}_w) \TbfM \partial^\tau \bfd^k_u.
\end{align*}
The first sum on the right-hand side is a telescoping sum in $k$, therefore, using \eqref{eq:qest}, we get
\begin{align*}
	\tau \sum_{k=q+1}^n (\bfde^k_w) \TbfM (\bfd^k_u - \eta \bfd^{k-1}_u) =&\ (1-\eta) \sum_{j=0}^{q-1} \sigma_j \big[(\bfe^{n-j}_w) \TbfM \bfd^n_u - (\bfe^{q-j}_w) \TbfM \bfd^q_u\big] \nonumber \\
	&\ - \tau \sum_{k=q+1}^n \sum_{j=0}^{q-1} \sigma_j (\bfe^{k-j-1}_w) \TbfM \partial^\tau \bfd^k_u \nonumber \\
	&\ + \tau \eta \sum_{k=q+1}^n \sum_{j=0}^{q-1} \sigma_j (\bfe^{k-j}_w) \TbfM \partial^\tau \bfd^k_u\nonumber \\
	\leq&\  c \| \bfd^n_u \|_{\ast,h}^2 + c \| \bfd^q_u \|_{\ast,h}^2 + \rho \sum_{k=n-q+1}^{n} \|\bfe^{k}_w \|_\bfK^2 + c \big( I_h^{q-1} + D_h^q \big) \nonumber \\
	&\ + \rho \tau \sum_{k=q+1}^n \| \bfe^k_w \|_\bfK^2 + c\tau \sum_{k=q+1}^n \| \partial^\tau \bfd^k_u \|_{\ast,h}^2 \nonumber \\
	\leq&\  \rho \sum_{\substack{k=n-q+1 \\ k\geq q+1}}^{n} \|\bfe^{k}_w \|_\bfK^2 + \rho \tau \sum_{k=q+1}^n \| \bfe^k_w \|_\bfK^2 \nonumber \\
	&\ + c\tau \sum_{k=q+1}^n \| \partial^\tau \bfd^k_u \|_{\ast,h}^2 
	+c\big(I_h^{q-1} + D_h^n\big) .
\end{align*}
The $\partial^\tau \bfd^k_u$ term is further estimated using Lemma~\ref{Lemma:partialtodot}, with the definition $\bfd^i_u=\bfd^u(i \tau)$ for $0 \leq i \leq q-1$, where $\bfd^u$ is the nodal vector of the defect in the continuous semi-discrete case. This yields the bound
\begin{align*}
	\tau \sum_{k=q+1}^n \|\partial^\tau \bfd^k_u \|_{\ast,h}^2 \leq c \tau \sum_{k=q}^n \|\mathbf{\dot d}^k_u \|_{\ast,h}^2 + c \tau  \sum_{i=1}^{q-1}, \|\partial^\tau \bfd^i_u \|_{\ast,h}^2 ,
\end{align*} 
which then gives \eqref{eq:boundproductrule}. 

Inserting the estimates \eqref{eq:ii.2.rhs1}, \eqref{eq:ii.2.rhs2}, and \eqref{eq:boundproductrule} into \eqref{eq:ii.2.lhs}, absorbing the term involving $\eta$ into the left-hand side,
using the estimates on the $G$-norm \eqref{eq:GnormEst}, and the estimate \eqref{eq:qest} for  $n=q$, we obtain the second part of the second energy estimate:
\begin{align}
	\sum_{k=n-q+1}^n \|\bfe_w^k\|_\bfA^2 + \tau \sum_{k=q+1}^n \|\bfde^k_u \|^2_\bfA \leq&\  c_0 \tau \sum_{k=q+1}^n \|\bfde^k_u \|_\bfM^2 + \widetilde \mu \tau \sum_{k=q+1}^n \| \bfde^k_u \|_\bfK^2 + \rho \sum_{\substack{k=n-q+1 \\ k\geq q+1\hphantom{n-}}}^{n} \|\bfe^{k}_w \|_\bfK^2 \nonumber \\
	&\ + \rho \tau \sum_{k=q+1}^{n} \| \bfe^k_w \|_\bfK^2+ \mu \tau \sum_{k=q+1}^{n-1} \| \bfe^k_u \|_\bfK^2 +c\big(I_h^{q-1} + D_h^n\big). \label{eq:ii.2bound}
\end{align}

The two parts of the second energy estimate are combined by
taking the linear combination of $2c_0$ times \eqref{eq:ii.1bound} and \eqref{eq:ii.2bound}. Absorptions and some rearrangement yields
\begin{align}
	\|  \bfe^n_w \|_\bfK^2 + \tau \sum_{k=q+1}^n \|\bfde^k_u \|_\bfK^2 \leq c \tau \sum_{k=q+1}^{n-1} \| \bfe^k_w \|_\bfK^2 + c \tau \sum_{k=q+1}^{n-1}  \| \bfe^k_u \|_\bfK^2 +c\big(I_h^{q-1} + D_h^n\big) . \label{eq:energyii}
\end{align}

\textit{(Combination)} We will now combine the two energy estimates obtained above, merging the two branches of Figure~\ref{fig:Scheme}. We add the first energy estimate \eqref{eq:energyi} and the second energy estimate \eqref{eq:energyii},
absorb the term involving $\bfde_u$, 
%
and use a discrete Gronwall inequality, together with the estimate \eqref{eq:qest} for $n=q$, to get the proposed energy estimate:
\begin{align*}
	\|  \bfe^n_u \|_\bfK^2 + \|  \bfe^n_w \|_\bfK^2 + \tau \sum_{k=q}^n \|\bfe^k_w \|_\bfK^2 + \tau \sum_{k=q}^n \|\bfde^k_u \|_\bfK^2 \leq c \big( I_h^{q-1} + D_h^n \big) .
\end{align*}

\textit{(Proving $t_{\max}=T$)}
It is only left to show that $t_{\max}=T$ for sufficiently small $h,\tau>0$. Let $n$ be the maximal natural number such that $(n-1) \tau \leq t_{\max}$ holds and $n \tau \leq T$. We then have, by the above inequality and the proposed bounds on the defects and initial values, that $\| \bfe^n_u \|_\bfK \leq C \, h^\kappa$. Then by the inverse estimate (see \cite[Theorem 4.5.11]{BrennerScott2008}), we have that
\begin{align*}
	\| \bfe^n_u \|_{L^\infty(\Om_h)} \leq&\  ch^{-\frac{d}{2}} \| \bfe^n_u \|_{L^2(\Om_h)} \leq ch^{-\frac{d}{2}} \| \bfe^n_u \|_\bfK \leq C \, h^{\kappa-\frac{d}{2}} \leq h^{\frac{\kappa}{2}-\frac{d}{4}}, 
\end{align*}
for sufficiently small $h$. Therefore $n\tau \leq t_{\max}$, which contradicts the maximality of $t_{\max}$, hence $t_{\max}=T$. \hfill
\end{proof}

\section{Consistency}
\label{section:Consistency}

The consistency analysis relies on an interplay of the geometric error estimates, the temporal error estimates and the error estimates for the spatial discretization, which are only recalled from \cite[Section 6]{HarderKovacs2022}. 
Some novel techniques are developed to estimate the $q$-step BDF time derivatives of the defects $\partial^\tau_q d^k_u$ and $\partial^\tau_q d^k_w$. 
In the following section we prove optimal-order bounds for the fully discrete defect terms. 

\subsection{Geometric error estimates} \label{Subsection:GeomErrorEst}
We recall, that the bounded domain $\Om \subset \mathbb{R}^d$ ($d=2,3$) has an at least $C^2$-boundary $\Ga$, and $\Ga_h=\partial \Om_h$, the boundary of the quasi-uniform triangulation $\Om_h$, is an interpolation of $\Ga$. By \cite[Theorem 5.9, Theorem 7.10]{ElliottRanner2021} (see also \cite{Dziuk1988}, \cite{Bernardi1989}, \cite{DziukElliott2013}, \cite{ElliottRanner2013}), we have that the finite element interpolation operator $\widetilde I_h v \in V_h$, with lift $I_h v =(\widetilde I_h v)^\ell \in V_h^\ell$, satisfies the following bounds:
\begin{lem}[interpolation errors] \label{Lemma:InterpolationError}
	For $ v \in H^2(\Om)$, such that $\gamma v \in H^2(\Ga)$, the piecewise linear finite element interpolation satisfies the following estimates, for $h\leq h_0$.\\
	(a) Interpolation error in the bulk:
	\begin{align*}
		\|v - I_h v \|_{L^2(\Om)} + h \| \nabla (v-I_hv) \|_{L^2(\Om)} \leq C \, h^2 \| v \|_{H^2(\Om)}.
	\end{align*}
	(b) Interpolation error on the surface (for $d=2,3$):
	\begin{align*}
		\| \gamma (v-I_h v) \|_{L^2(\Ga)} + h \| \nabla_{\Ga} (v-I_h v) \|_{L^2(\Ga)} \leq C \, h^2 \| \gamma v \|_{H^2(\Ga)}.
	\end{align*}
\end{lem}
For the surface estimate we note here, that $\Ga$ is a $d-1 \leq 3$ dimensional surface in $\R^d$, whence (b) does not hold in higher dimensions, cf., for instance, \cite[Lemma~4.3]{DziukElliott2013}. 

The following error estimates for the Ritz map \eqref{eq:DefRitzmap} are from \cite[Lemma~3.1, 3.3]{KovacsLubich2017}. We note that compared to the Ritz map in \cite{KovacsLubich2017} \eqref{eq:DefRitzmap} contains some lower order bulk terms.
\begin{lem}[Ritz map errors] \label{Lemma:RitzMapError}
	For any $v \in H^2(\Om)$, with $\gamma v \in H^2(\Ga)$, the error of the Ritz map \eqref{eq:DefRitzmap} satisfies the following bounds, for $h \leq h_0$,
	\begin{align*}
		\| v- R_hv\|_{H^1(\Om)} + \| \gamma (v-R_hv) \|_{H^1(\Ga)} \leq C \, h \Big(\| v \|_{H^2(\Om)} + \| \gamma v \|_{H^2(\Ga)}\Big), \\
		\| v- R_hv\|_{L^2(\Om)} + \| \gamma (v-R_hv) \|_{L^2(\Ga)} \leq C \, h^2 \Big(\| v \|_{H^2(\Om)} + \| \gamma v \|_{H^2(\Ga)}\Big),
	\end{align*}
	where the constant $C > 0$ is independent of $h$.
\end{lem}
The last result, from \cite[Lemma 6.2]{KovacsLubich2017} (combined with \cite[Lemma 6.3]{ElliottRanner2013}), estimates the error for the bilinear forms $a,m$ and $a_h,m_h$.
\begin{lem}[geometric approximation errors] \label{Lemma:GeomAppError}
	The bilinear forms $a,m$ and their discrete counterparts $a_h,m_h$, satisfy the following estimates, for $h \leq h_0$ and for any $v_h,w_h \in V_h$,
	\begin{align*}
		|a(v^\ell_h,w^\ell_h) - a_h(v_h,w_h)| \leq&\  C \, h \| \nabla v^\ell_h \|_{L^2(\Om)} \| \nabla w^\ell_h \|_{L^2(\Om)} \nonumber \\
		&\ + C \, h^2 \Big(\| \nabla v_h^\ell \|_{L^2(\Om)} \| \nabla w_h^\ell \|_{L^2(\Om)} + \| \nabla_\Ga v_h^\ell \|_{L^2(\Ga)} \| \nabla_\Ga w_h^\ell \|_{L^2(\Ga)}\Big), \\
		|m(v^\ell_h,w^\ell_h) - m_h(v_h,w_h)| \leq&\  C \, h^2 \|  v^\ell_h \|_{H^1(\Om)} \|  w^\ell_h \|_{H^1(\Om)} \nonumber \\
		&\ + C \, h^2 \Big(\|  v_h^\ell \|_{L^2(\Om)} \| w_h^\ell \|_{L^2(\Om)} + \| \gamma v_h^\ell \|_{L^2(\Ga)} \| \gamma w_h^\ell \|_{L^2(\Ga)}\Big).
	\end{align*}
	A combination of these estimates yields a similar estimate for the bilinear forms $a^\ast$ and $a^\ast_h$.
\end{lem}

\begin{rem} \label{rem:normequivalence}
As a consequence we also have the $h$-uniform equivalence of the norms $\| \cdot \|$ and $\| \cdot \|_h$ on $V_h^\ell$ and $V_h$, and of the norms $| \cdot |$ and $| \cdot |_h$ on $V_h^\ell$ and $V_h$.
\end{rem}

\subsection{Temporal error estimates}
We recall the definition of the discrete derivative and the extrapolation,
with the coefficients \eqref{eq:BdfCoeff}:  
\begin{align*}
	\partial^\tau_q u^n = \frac{1}{\tau} \sum_{j=0}^q \delta_j u^{n-j} ,\qquad	\widetilde u^n = \sum_{j=0}^{q-1} \gamma_j u^{n-1-j} .
\end{align*}
We can define this in a more general way for a function $y:[0,T] \to \mathbb{R}$ as 
\begin{align}
	\partial^\tau_q y(s)= \frac{1}{\tau} \sum_{i=0}^q \delta_i y(s-i\tau), \qquad \widetilde y (s) = \sum_{j=0}^{q-1} \gamma_j y(s-j\tau-\tau),\label{eq:contdefpartialtilde}
\end{align}
for $s\geq q\tau$, which is consistent with the previous definition for $s=n\tau$ and the natural choice $u^j=y(j\tau)$ for $j=n-q,\dotsc,n$. 

The main tool in analyzing the temporal error is the following result, based on the Plano kernel representation of the error for a multi step method, see, e.g., \cite[Section~III.2]{HairerNorsettWanner1993}.
\begin{lem}
\label{Lemma:PeanoKernel}
	(i) Let $y\colon[0,T] \to \mathbb{R}$ be $p$ times continuously differentiable, with $y^{(p)}$ absolutely continuous, for $1 \leq p \leq q \leq 5$, and let $t^\ast \in [q\tau,T]$, then 
	\begin{align*}
		\partial^\tau_q y(t^\ast) - \frac{\d}{\d t}y(t^\ast)=\tau^p \int_{0}^{q} K_q(s) y^{(p+1)}(t^\ast-s\tau) ~\d s ,
	\end{align*}
	where
	\begin{align*}
		K_{q}(s)= - \frac{1}{p!} \sum_{i=0}^q \delta_i (s-i)_-^p , \qquad \text{ with } \qquad (x)_- := \begin{cases}
			x, & x \leq 0 , \\
			0, & x \geq 0 .
		\end{cases}
	\end{align*}
	(ii) Let $y\colon[0,T] \to \mathbb{R}$ be $p-1$ times continuously differentiable, with $y^{(p-1)}$ absolutely continuous, for $1 \leq p \leq q \leq 5$, and let $t^\ast \in [q\tau,T]$, then 
	\begin{align*}
		\widetilde y(t^\ast) - y(t^\ast)=\tau^p \int_{0}^{q} \widetilde K_q(s) y^{(p)}(t^\ast-s\tau) ~\d s,
	\end{align*}
	where 
	\begin{align*}
		\widetilde K_q(s) = \sum_{i=0}^{q-1} \gamma_i \frac{(s-(i+1))_-^{p-1}}{(p-1)!}.
	\end{align*}
\end{lem}
\begin{proof}
	The first equation is obtained by using Taylor's theorem, see, e.g., \cite[Section~III.2]{HairerNorsettWanner1993}. The second equation follows by a similar argument. \hfill
\end{proof}

Since as part of the consistency analysis we have to bound the discrete derivative of the defect terms $d_u$ and $d_w$, we will use Lemma~\ref{Lemma:PeanoKernel} in the sub-optimal case $p=1$, for both the error in the discrete derivative $y=\partial^\tau_q u_\ast - \frac{\d}{\d t} u_\ast$ and the error in the extrapolation $y=\widetilde u_\ast - u_\ast$. 

Therefore, we have to extend the definition \eqref{eq:contdefpartialtilde} of $\partial^\tau_q u_\ast(t)$ to $t<q\tau$. It is natural to define $\partial^\tau_q u_\ast(t_i)= \frac{\d}{\d t} u_\ast(t_i)$, for $0 \leq i \leq q-1$, and since we need the extension to be in $C^{2}$, we define
\begin{align}
	\partial^\tau_q u_\ast(t):=\begin{cases}
		\frac{1}{\tau} \sum_{j=0}^{q} \delta_j u_\ast(t-j\tau)  , ~& \text{for} \, q\tau \leq t \leq T \\
		p(t)  , ~& \text{for} \, (q-1)\tau \leq t \leq q\tau \\
		\frac{\d}{\d t} u_\ast(t)  , ~& \text{for} \, 0 \leq t \leq (q-1)\tau ,
	\end{cases}, \label{eq:GenDefPartial}
\end{align}
where $ p\colon [(q-1)\tau,q\tau] \to \mathbb{R}$ is defined as a \emph{quintic Hermite interpolation}
, then we have $\partial^\tau_q u_\ast\in C^{2}([0,T])$ and obtain
\begin{align}
	\Big\| \int_{q-1}^q K_q(k-s)\frac{\d^2}{\d t^2} \Big( p- \frac{\d}{\d t} u_\ast \Big)(s\tau) ~\d s \Big\|_{\ast,h}   \leq C \, \tau^{q-1}. \label{eq:HermiteCondPartial}
\end{align}

We construct $p:[(q-1)\tau,q\tau] \to \R$ as the quintic Hermite polynomial, corresponding to the nodes 
\begin{equation*}
	\bigg(x_0 := (q-1) \tau ~ , ~ y_0^j := \frac{\d^{j+1}}{\d t^{j+1}} u_\ast((q-1)\tau)\bigg) \quad \text{ and } \quad \bigg(x_1 := q \tau ~ , ~ y_1^j := \frac{\d^j}{\d t^j} \frac{1}{\tau} \sum_{j=0}^{q} \delta_j u_\ast(q\tau-j\tau)\bigg),
\end{equation*}
for $j=0,1,2$. 

The standard quintic Hermite polynomial $\widetilde p : [(q-1)\tau,q\tau] \to \R$ approximating $\frac{\d}{\d t} u_\ast$ satisfies by standard error analysis (see, e.g., \cite[Theorem 3.1]{WongAgarwal1989}), for $u_\ast \in C^{q+1}([(q-1)\tau,q\tau])$ and with $\varphi_h \in V_h$, the estimate
\begin{align*}
	m_h\Big(\int_{q-1}^q K_q(k-s) \frac{\d^2}{\d t^2} \Big(\widetilde p - \frac{\d}{\d t} u_\ast\Big),\varphi_h\Big) \leq C \, \Big\| \frac{\d^2}{\d t^2} \Big(\widetilde p - \frac{\d}{\d t} u_\ast\Big) \Big\|_{L^\infty(\Omega_h)} \| \varphi_h \|_h \leq C \, \tau^{q-1} \|\varphi_h \|_h
\end{align*}
We compare the standard Hermite polynomial $\widetilde p$ to $p$ by the representation in Hermite basis functions and their bounds from \cite[Lemma 2.1]{WongAgarwal1989}:
\begin{align*}
	m_h\Big(\int_{q-1}^{q} K_q(k-s) \frac{\d^2}{\d t^2} (p-\widetilde p)(s\tau) \d s,\varphi_h\Big)\leq&\ C ~ m_h\Big(\int_{q-1}^{q} \frac{\d^2}{\d t^2} (p-\widetilde p)(s\tau) \d s,\varphi_h\Big) \\
	\leq&\ C \max_{j=0,1,2} \tau^{j-1} m_h\Big(\frac{1}{\tau}\sum_{j=0}^q \delta_j u_\ast^{(j)}((q-j)\tau)  - u_\ast^{(j+1)}(q\tau), \varphi_h\Big) \\
	\leq&\ C \max_{j=0,1,2} \tau^{q-1} m_h(\int_0^q \widetilde K_q(s) u_\ast^{(q+1)}(q\tau-s\tau)\d s, \varphi_h)\\
	\leq&\ C \, \tau^{q-1} \| \varphi_h \|_h , 
\end{align*}
where we again used Lemma~\ref{Lemma:PeanoKernel}, and in the last inequality the estimate \eqref{eq:Rubound} to control the term in the Ritz map of the exact solution by the regularity assumptions \eqref{eq:RegularityAss}, via the Ritz map error estimates of Lemma~\ref{Lemma:RitzMapError}. Combining the last two estimates we arrive at the desired estimate \eqref{eq:HermiteCondPartial}.

Similarly we need to extend the definition \eqref{eq:contdefpartialtilde} for the extrapolation, which for $0 \leq t \leq (q-1)\tau$ is the obvious choice $u_\ast(t)$:
\begin{align}
	\widetilde u_\ast(t) = \begin{cases}
		\sum_{j=1}^{q} \gamma_j u_\ast(t-j\tau)  , ~& \text{for} \, q\tau \leq t \leq T , \\
		r(t)  , ~& \text{for} \, (q-1)\tau \leq t \leq q\tau , \\
		u_\ast(t)  , ~& \text{for} \, 0\leq t \leq (q-1) \tau ,
	\end{cases} \label{eq:defExtrapol}
\end{align}
where we choose $r \colon [(q-1)\tau, q\tau] \to \mathbb{R}$ again as the quintic Hermite polynomial, such that $\widetilde u_\ast \in C^{2}([0,T])$ and 
\begin{align}
	\Big\| \int_{q-1}^q K_q(k-s) \frac{\d^2}{\d t^2}  \big(r - u_\ast \big)(s\tau) ~\d s\Big\|_{\ast,h}   \leq C \, \tau^{q-1}. \label{eq:HermiteCondtilde}
\end{align}

\subsection{Fully discrete defect bounds}
With this preparation we can start to estimate the different error defect terms, which appear in the Stability estimates.
\begin{prop} \label{Prop:ConsistencyEst}
	Let $u$ and $w$ solve \eqref{eq:strong2CH} and satisfy the regularity condition \eqref{eq:RegularityAss}. Then we have, for $q \leq k \leq n$ and $1\leq i \leq q-1$, the following defect estimates:
	\begin{alignat}{3}
		\| d^k_u \|_{\ast,h} \leq&\  C \, (h^2 + \tau^q),\qquad &&\| \partial^\tau_q d^k_u \|_{\ast,h} \leq&&\  C \, (h^2 + \tau^q), \label{eq:du}\\
		\| d^k_w \|_{\ast,h} \leq&\  C \, (h^2 + \tau^q), \qquad &&\| \partial^\tau_q d^k_w \|_{\ast,h} \leq&&\  C \, (h^2 + \tau^q), \label{eq:dw}\\
		\|  d^{q-1}_w \|_{\ast,h} \leq&\  C \, h^2 , \qquad &&\| \partial^\tau d^i_w \|_{\ast,h} \leq&&\  C \, h^2 .\label{eq:pdi}
	\end{alignat}
	Therefore, there holds
	\begin{align*}
		D_h^n\leq C \, (h^2 + \tau^{q})^2.
	\end{align*}
\end{prop}
\begin{proof} 
	The defect $d_u^k$ is defined by the equation \eqref{eq:dCH} in the abstract form using bilinear forms. By adding intermediate terms, such as the original equation for the exact solution, we split this into a temporal and a spatial defect:
	\begin{align*}
		m_h(d_u^k,\varphi_u)=&\ m_h(\partial^\tau_q u_\ast^k,\varphi_u) + a_h(w_\ast^k,\varphi_u)- m(\dot u(t_k),\varphi_u^\ell) - a(w(t_k),\varphi_u^\ell) \\
		=&\ m_h(\partial^\tau_q u_\ast^k,\varphi_u) - m_h\Big(\frac{\d }{\d t} u_\ast^k,\varphi_u\Big) + m_h(d_u(t_k),\varphi_u),
	\end{align*}	
	for $t_k=k\tau$. Analogously, we have for the defect $d_w^k$:
	\begin{align*}
		m_h(d_w^k,\varphi_w)=&\ m_h(W'(u_\ast(t_k)) - W'(\widetilde u_\ast (t_k)), \varphi_w) + m_h(d_w(t_k).
	\end{align*}
	Recalling the optimal-order defect bounds for the spatial defect from \cite[Proposition 6.1]{HarderKovacs2022} we have, for $0 \leq t \leq T$,
	\begin{equation}\label{eq:defectcont}
		\begin{aligned}
			&\| d_u(t) \|_{\ast,h} \leq C \, h^2 , \qquad &&\| \dot d_u(t) \|_{\ast,h} \leq C \, h^2,\qquad \text{and } &&\| \ddot d_u(t) \|_{\ast,h} \leq C \, h^2, \\
			&\| d_w(t) \|_{\ast,h} \leq C \, h^2 , \qquad &&\| \dot d_w(t) \|_{\ast,h} \leq C \, h^2,\qquad \text{and } &&\| \ddot d_w(t) \|_{\ast,h} \leq C \, h^2,
		\end{aligned}
	\end{equation}
	The estimates for the second derivative follow exactly in the same way, but need the exact solution to be one derivative more regular.
	
	We start with the last two estimates \eqref{eq:pdi}. Using the fact that, for $0 \leq i \leq q-1$, we defined $d^i_w=d_w(i \tau)$, we have by \eqref{eq:defectcont}, and the mean value theorem
	\begin{align*}
		&\| d^{q-1}_u \|_{\ast,h} = \| d_w((q-1)\tau)\|_{\ast,h} \leq C \, h^2 , \\
		&\| \partial^\tau d^i_w \|_{\ast,h} = \| \dot d_w(\zeta) \|_{\ast,h} \leq C \, h^2 \qquad \text{ for some } \ \zeta \in [(i-1) \tau, i \tau].
	\end{align*} 
	
	We turn to \eqref{eq:du}, by using Lemma~\ref{Lemma:PeanoKernel}, the spatial defect bounds \eqref{eq:defectcont}, and the fact that $\frac{\d}{\d t} (\widetilde{R}_h u) = \widetilde{R}_h \dot u$,
	\begin{align}
		m_h(d^k_u,\varphi_u)=&\  m_h\Big(\partial^\tau_q u_\ast(t_k)-\frac{\d}{\d t} u_\ast(t_k),\varphi_u\Big) + m_h(d_u(t_k),\varphi_u) \nonumber \\
		=&\  m_h\Big( \tau^q \int_{0}^{q} K_q(s) u_\ast^{(q+1)}(t_k-s\tau) ~\d s ,\varphi_u\Big) + m_h(d_u(t_k),\varphi_u) \nonumber \\
		\leq&\  C \, \tau^q |m_h(\widetilde{R}_h u^{(q+1)},\varphi_u) |_{L^\infty} + C \, h^2 \| \varphi_u \|_h, \label{eq:du1}
	\end{align}
	we further estimate by the geometric error estimates of Section~\ref{Subsection:GeomErrorEst}, for $0 \leq t \leq T$ and $0 \leq j \leq q+2$,
	\begin{align}
		|m_h(\widetilde{R}_h u^{(j)}(t),\varphi_u)|\leq&\  |m_h(\widetilde{R}_h u^{(j)}(t),\varphi_u) - m( R_h u^{(j)}(t),\varphi_u^\ell)| + |m( R_h u^{(j)}(t) - u^{(j)}(t),\varphi_u^\ell)|\nonumber \\
		&\ + |m(u^{(j)}(t),\varphi_u^\ell)| \nonumber \\
		\leq&\  C \, h \|  R_h u^{(j)} (t) \|_{L^2(B^\ell_h)} \| \varphi_u^\ell \|_{L^2(B^\ell_h)} + C \, h^2 |  R_h u^{(j)} (t) | | \varphi_u^\ell | \nonumber \\
		&\  + C \, h^2 \big(\| u^{(j)}(t)\|_{H^2(\Om)} + \| \gamma u^{(j)}(t) \|_{H^2(\Ga)}\big) | \varphi_u^\ell | + | u^{(j)}(t) | |\varphi_u^\ell| \nonumber \\
		\leq&\  ch^2 \| R_h u^{(j)}(t) - u^{(j)}(t)\| \| \varphi_u \|_h + ch^2 \| u^{(j)}(t)\| \| \varphi_u \|_h + C \, \| \varphi_u \|_h \nonumber \\
		\leq&\  ch^2 (1+ch) \big(\| u^{(j)}(t) \|_{H^2(\Om)} + \| \gamma u^{(j)} (t)\|_{H^2(\Ga)}\big) \| \varphi_u \|_h + C \, \| \varphi_u \|_h \nonumber \\
		\leq&\  C \, \| \varphi_u \|_h, \label{eq:Rubound}
	\end{align}
	since the solution $u$ is sufficiently regular, and where the constant only depends on $u$. Combining \eqref{eq:du1} and \eqref{eq:Rubound}, yields the first estimate in \eqref{eq:du}. 
	
	In order to analyze the discrete derivative of $d^k_u$, we consider
	\begin{align}
		m_h(\partial^\tau_q d^k_u, \varphi_u)=&\  m_h\bigg(\partial^\tau_q \Big( \partial^\tau_q u_\ast(t_k)\Big) - \partial^\tau_q \Big( \frac{\d}{\d t} u_\ast(t_k)\Big), \varphi_u\bigg) + m_h(\partial^\tau_q d_u(t_k),\varphi_u) \nonumber \\
		=&\  m_h\bigg(\Big(\partial^\tau_q - \frac{\d}{\d t}\Big) \Big( \partial^\tau_q u_\ast(t_k) - \frac{\d}{\d t} u_\ast(t_k)\Big), \varphi_u\bigg) + m_h\bigg(\frac{\d}{\d t}\Big(\partial^\tau_q u_\ast(t_k) - \frac{\d}{\d t} u_\ast(t_k)\Big),\varphi_u\bigg) \nonumber \\
		&\ + m_h\Big(\partial^\tau_q d_u(t_k) - \frac{\d}{\d t} d_u(t_k),\varphi_u\Big) + m_h( \dot d_u(t_k),\varphi_u) \nonumber \\
		=: &\ I + II + III + IV. \label{eq:partialduest}
	\end{align}
	We estimate these terms, by using Lemma~\ref{Lemma:PeanoKernel} (for $p=q$ and $p=1$), estimate \eqref{eq:Rubound} and the estimates on the spatial defects \eqref{eq:defectcont}:
	\begin{subequations}
	\label{eq:partialII-IV}
		\begin{align}
			\label{eq:estddu}
			I= &\ m_h\bigg(\Big(\partial^\tau_q - \frac{\d}{\d t}\Big) \Big( \partial^\tau_q u_\ast(t_k) - \frac{\d}{\d t} u_\ast(t_k) \Big), \varphi_u\bigg) \leq C \, \tau^q \| \varphi_u \|_h, \\
			II = &\ m_h\bigg(\frac{\d}{\d t}\Big(\partial^\tau_q u_\ast(t_k) - \frac{\d}{\d t} u_\ast(t_k)\Big),\varphi_u\bigg) = m_h\Big( \tau^q \int_{0}^{q} K_q(s) u_\ast^{(q+2)}(t_k-s\tau) ~\d s ,\varphi_u\Big) \leq C \, \tau^q \| \varphi_u \|_h, \\
			III = &\ m_h\Big(\partial^\tau_q d_u(t_k) - \frac{\d}{\d t} d_u(t_k),\varphi_u\Big) = m_h\Big( \tau^1 \int_{0}^{q} K_q(s) d_u^{(2)}(t_k-s\tau) ~\d s ,\varphi_u\Big) \leq C \, h^2 \| \varphi_u \|_h, \\
			IV = &\ m_h( \dot d_u(t_k),\varphi_u)\leq C \, h^2 \| \varphi_u \|_h. 
		\end{align}
	\end{subequations}
	In order to establish \eqref{eq:estddu}, we use Lemma~\ref{Lemma:PeanoKernel}, with $y=\partial^\tau_q u_\ast (t_k)- \frac{\d}{\d t} u_\ast(t_k)$ and the sub-optimal estimate $p=1$: 
	\begin{align}
		m_h\bigg(\Big(\partial^\tau_q - \frac{\d}{\d t}\Big) \Big( \partial^\tau_q u_\ast(t_k) - \frac{\d}{\d t} u_\ast(t_k)\Big), \varphi_u\bigg) =&\ m_h\Big( \tau^1 \int_{0}^{q} K_q(s) y^{(2)}(t_k-s\tau) ~\d s ,\varphi_u\Big) \label{eq:estddu1}
	\end{align}
	Then we estimate with Lemma~\ref{Lemma:PeanoKernel}, for $p=q-1$, and the estimate \eqref{eq:Rubound}, for the interval $t_{q}\leq t \leq T$:
	\begin{align}
		m_h\bigg( \frac{\d^2}{\d t^2}  \Big(\partial^\tau_q - \frac{\d}{\d t}\Big) u_\ast(t),\varphi_u\bigg)=&\  \tau^{q-1} m_h \bigg(  \frac{\d^2}{\d t^2}  \int_0^q K_q(s) u_\ast^{(q)}(t-s\tau) ~\d s , \varphi_u\bigg) \nonumber \\
		\leq&\  C \, \tau^{q-1} |m_h(\widetilde{R}_h u^{(q+2)},\varphi_u)|_{L^\infty([0,T])} \nonumber \\
		\leq&\  C \, \tau^{q-1} \| \varphi_u \|_h. \label{eq:estddu2}
	\end{align}
	which together with \eqref{eq:HermiteCondPartial}, for the interval $(q-1)\tau \leq t < q\tau$ and by the definition \eqref{eq:GenDefPartial}, for the interval $0\leq t <(q-1)\tau $, results in
	\begin{align*}
		m_h\bigg(\int_0^q K_q(s) \frac{\d^2}{\d t^2}  \Big(\partial^\tau_q u_\ast(t_k-s\tau) - \frac{\d}{\d t} u_\ast(t_k-s\tau)\Big) ~\d s, \varphi_u\bigg) \leq C \, \tau^{q-1} \| \varphi_u \|_h.
	\end{align*}
	Inserting this in \eqref{eq:estddu1} we arrive at \eqref{eq:estddu}.
	
	Together with \eqref{eq:partialduest}, the estimates of \eqref{eq:partialII-IV} yield the second estimate in \eqref{eq:du}.
	
	The estimates for $d^k_w$ and $\partial^\tau_q d^k_w$ follow in an analogous ways, except for the non-linear term. For which we have 
	\begin{align*}
		m_h\Big(W'\big(u_\ast(t_k)\big) - W'\big(\widetilde u_\ast (t_k)\big), \varphi_w\Big) \leq  C \, m_h(u_\ast(t_k) - \widetilde u_\ast(t_k), \varphi_w), 
	\end{align*}
	by the local Lipschitz continuity of $W'$. Then, similarly as above, we use Lemma~\ref{Lemma:PeanoKernel} and the estimate \eqref{eq:Rubound}, to obtain
	\begin{align}
		m_h(u_\ast(t_k) - \widetilde u_\ast(t_k), \varphi_w) =&\  \tau^q m_h\Big(\int_0^q \widetilde K_q(s) u_\ast^{(q)}(t_k-s) ~\d s, \varphi_w\Big) \nonumber \\
		\leq&\  C \, \tau^q | m_h(\widetilde{R}_h u^{(q)}, \varphi_w) |_{L^\infty([0,T])} \nonumber \\
		\leq&\  C \, \tau^q \| \varphi_w \|_h, \label{eq:tildeu_ast}
	\end{align}
	and together with the estimates of the spatial defect \eqref{eq:defectcont}, we obtain the first estimate \eqref{eq:dw}. 
	
	For the discrete derivative of the expression we consider
	\begin{align}
		m_h(\partial^\tau_q d^k_w, \varphi_w)  =&\  m_h\Big(\partial^\tau_q \big(W'\big(\widetilde u_\ast(t_k)\big)  - W'\big(u_\ast(t_k)\big)\big), \varphi_w\Big)  \nonumber \\
		&\ + \bigg( m_h\bigg(\Big(\partial^\tau_q - \frac{\d}{\d t}\Big) d_w(t_k), \varphi_w\bigg) + m_h(\dot d_w(t_k),\varphi_w) \bigg) \nonumber \\
		=:&\ I + II. \label{eq:ddw1}
	\end{align}
	By analogous estimates as in \eqref{eq:partialII-IV}, we have
	\begin{align}
		II = m_h\bigg(\Big(\partial^\tau_q - \frac{\d}{\d t}\Big) d_w(t_k), \varphi_w\bigg) + m_h(\dot d_w(t_k),\varphi_w) \leq C \, h^2 \| \varphi_w \|_h. \label{eq:partialdw1}
	\end{align}
	Therefore, it is left to estimate the term $I$ in \eqref{eq:ddw1}. We have that
	\begin{align}
		\partial^\tau_q \big(W'(\widetilde u_\ast(t_k))  - W'(u_\ast(t_k))\big) = \sum_{j=k-q+1}^k \sigma_{k-j} \partial^\tau \big(W'(\widetilde u_\ast(t_j))  - W'(u_\ast(t_j))\big), \label{eq:ddw2}
	\end{align}
	and by setting $u^\ast_\theta=u_\ast + \theta (\widetilde u_\ast - u_\ast)$, for $1 \leq i \leq k$, we have
	\begin{align*}
		\partial^\tau \big(W'\big(\widetilde u_\ast(t_i)\big)  - W'\big(u_\ast(t_i)\big)\big) =&\  \int_0^1 \frac{d}{d\theta} \partial^\tau \Big(W'\big(u^\ast_\theta (t_i)\big)\Big) ~\d \theta \nonumber \\
		=&\  \int_0^1 \partial^\tau  \Big(W''\big(u^\ast_\theta (t_i)\big) \big(\widetilde u_\ast(t_i) - u_\ast(t_i)\big)\Big) ~\d \theta \nonumber \\
		=&\  \int_0^1 W''\big(u^\ast_\theta (t_i)\big) \partial^\tau \big(\widetilde u_\ast(t_i) - u_\ast(t_i)\big) ~\d \theta \nonumber \\
		&\ + \int_0^1 \partial^\tau W''\big(u^\ast_\theta (t_i)\big) \big(\widetilde u_\ast(t_{i-1}) - u_\ast(t_{i-1})\big) ~\d \theta \nonumber \\
		\leq&\ C | \partial^\tau (\widetilde u_\ast(t_i) - u_\ast(t_i))| + C |\widetilde u_\ast(t_{i-1}) - u_\ast(t_{i-1})| ,
	\end{align*}
	where we used the local Lipschitz continuity (as in the stability proof for a large enough compact neighborhood around the exact solution, see \eqref{eq:r1.1}), together with the $L^\infty$-bounds
	\begin{align*}
		&\| \partial^\tau u_\ast  \|_{L^\infty([\tau,k\tau],L^\infty(\Om_h))} \leq  \| \dot u_\ast  \|_{L^\infty([0,T],L^\infty(\Om_h))} \leq C ,\\
		&\| \widetilde u_\ast - u_\ast \|_{L^\infty([0,(k-1)\tau],L^\infty(\Om_h))} \leq  C ,
	\end{align*}
	which follow by Lemma~\ref{Lemma:PeanoKernel} and the regularity assumptions \eqref{eq:RegularityAss}. Therefore, we obtain
	\begin{align}
		m_h\Big(\partial^\tau \big(W'\big(\widetilde u_\ast(t_i)\big)  - W'\big(u_\ast(t_i)\big)\big),\varphi_w\Big)\leq&\  C \, \| \varphi_w \|_h \Big( \| \partial^\tau (\widetilde u_\ast(t_i) - u_\ast(t_i))\|_{\ast,h} + \|\widetilde u_\ast(t_{i-1}) - u_\ast(t_{i-1})\|_{\ast,h}\Big). \label{eq:ddw3}
	\end{align}
	Since $\widetilde u_\ast(t_i) - u_\ast(t_i)=0$, for $0 \leq i \leq q-1$, we combine \eqref{eq:ddw2} and \eqref{eq:ddw3}, and use Lemma~\ref{Lemma:partialtodot}, to get
	\begin{align*}
		m_h\Big(\partial^\tau_q \big(W'\big(\widetilde u_\ast(t_i)\big)  - W'\big(u_\ast(t_i)\big)\big),\varphi_w\Big)\leq&\  C \, \| \varphi_w \|_h \Big(  \sum_{\substack{j=k-q+1 \\ j\geq q\hphantom{k-+1}}}^{k} \| \partial^\tau_q (\widetilde u_\ast(t_j) - u_\ast(t_j))\|_{\ast,h} + \|\widetilde u_\ast(t_{j}) - u_\ast(t_{j})\|_{\ast,h} \Big). 
	\end{align*}
	Then by estimate \eqref{eq:tildeu_ast}, we have
	\begin{align*}
		m_h\Big(\partial^\tau_q \big(W'\big(\widetilde u_\ast(t_i)\big)  - W'\big(u_\ast(t_i)\big)\big),\varphi_w\Big)\leq&\  C \, \| \varphi_w \|_h \sum_{\substack{j=k-q+1 \\ j\geq q\hphantom{k-+1}}}^{k} \| \partial^\tau_q (\widetilde u_\ast(t_j) - u_\ast(t_j))\|_{\ast,h} + c \tau^q \| \varphi_w \|_h.
	\end{align*}
	And therefore, we further estimate, for $q\leq k \leq n$,
	\begin{align}
		m_h\Big(\partial^\tau_q \big(\widetilde u_\ast(t_k)  - u_\ast(t_k)\big), \varphi\Big) = m_h\bigg(\Big(\partial^\tau_q - \frac{\d}{\d t}\Big) \Big(\widetilde u_\ast(t_k) - u_\ast(t_k)\Big), \varphi\bigg) + m_h\bigg(\frac{\d}{\d t} \Big(\widetilde u_\ast(t_k) - u_\ast(t_k)\Big), \varphi\bigg) . \label{eq:partialdw2}
	\end{align}
	The second term we can estimate, using Lemma~\ref{Lemma:PeanoKernel} and \eqref{eq:Rubound}, by
	\begin{align}
		m_h\bigg(\frac{\d}{\d t} \Big(\widetilde u_\ast(t_k) - u_\ast(t_k)\Big), \varphi\bigg) =&\  \tau^q m_h\Big(\int_0^q \widetilde K_q(s) u_\ast^{(q+1)}(t_k - s\tau) ~\d s, \varphi\Big) \nonumber \\
		\leq&\  C \, \tau^q | m_h(\widetilde{R}_h u^{(q+1)}, \varphi ) |_{L^\infty([0,T])} \nonumber \\
		\leq&\  C \, \tau^q \| \varphi \|_h. \label{eq:partialdw3}
	\end{align}
	
	Similar to the estimate for $\partial^\tau_q d_u^k$, for the first term we use Lemma~\ref{Lemma:PeanoKernel}, for the function $y=\widetilde u_\ast(t) - u_\ast(t)$ in the sub-optimal case $p=1$, and obtain
	\begin{align}
		m_h\bigg(\Big(\partial^\tau_q - \frac{\d}{\d t}\Big) \Big(\widetilde u_\ast(t_k) - u_\ast(t_k)\Big), \varphi\bigg)=&\ m_h\Big( \tau^1 \int_{0}^{q} K_q(s) y^{(2)}(t_k-s\tau) ~\d s ,\varphi_u\Big). \label{eq:almostestddw}
	\end{align}
	For $q\tau \leq t \leq T$, we have with Lemma~\ref{Lemma:PeanoKernel}, for $p=q-1$, and the estimate \eqref{eq:Rubound}, which yield 
	\begin{align}
		m_h\bigg(  \frac{\d^2}{\d t^2}  \Big(\widetilde u_\ast(t) - u_\ast(t)\Big), \varphi\bigg) =&\  \tau^{q-1} m_h\bigg(  \frac{\d^2}{\d t^2}  \int_0^q \widetilde K_q(s) u_\ast^{(q-1)}(t-s\tau) ~\d s,\varphi\bigg) \nonumber \\
		\leq&\  C \, \tau^{q-1} | m_h(\widetilde{R}_h u^{(q+1)},\varphi) |_{L^\infty([0,T])} \nonumber \\
		\leq&\  C \, \tau^{q-1} \| \varphi \|_h. \label{eq:ddw4}
	\end{align}
	Combining \eqref{eq:ddw4} with the estimate \eqref{eq:HermiteCondtilde}, for the interval $(q-1)\tau \leq t < q \tau$, using the fact that $\widetilde u_\ast = u_\ast$ in the interval $0\leq t \leq (q-1)\tau$ by the definition \eqref{eq:defExtrapol}, and inserting into \eqref{eq:almostestddw}, we arrive at
	\begin{align}
		m_h\bigg(\Big(\partial^\tau_q - \frac{\d}{\d t}\Big) \Big(\widetilde u_\ast(t_k) - u_\ast(t_k)\Big), \varphi\bigg) \leq&\  C \, \tau^{q} \| \varphi \|_h. \label{eq:partialdw4}
	\end{align}
	Combining all the above estimates \eqref{eq:partialdw1}, \eqref{eq:partialdw2}, \eqref{eq:partialdw3}, and \eqref{eq:partialdw4} we arrive at
	\begin{align*}
		m_h(\partial^\tau_q d^k_w, \varphi_w)\leq&\  C (\tau^q+h^2) \| \varphi_w \|_h,
	\end{align*} 
	which yields the last estimate in \eqref{eq:dw}. \hfill
\end{proof}

\section{Optimal-order fully discrete error estimates}
\label{section:MainProof}

The last two sections were devoted to establishing stability and consistency for the Cahn--Hilliard equation with Cahn--Hilliard-type dynamic boundary conditions. We now combine the results to prove our main result Theorem~\ref{Thm:ConvergenceMainResult}.

\begin{proof}[Proof of Theorem~\ref{Thm:ConvergenceMainResult}]
	The error between the lifted numerical solutions at time $t_n$ and the exact solutions decomposes to
	\begin{equation} \label{eq:MainProofdecomp}
		\begin{aligned}
			(u^n_h)^\ell - u(t_n)=&\  (u^n_h - \widetilde{R}_h u(t_n))^\ell + (R_h u(t_n) - u(t_n)) =  (e^n_u)^\ell + (R_h u(t_n) - u(t_n)), \\
			(w^n_h)^\ell - w(t_n)=&\  (w^n_h - \widetilde{R}_h w(t_n))^\ell + (R_h w(t_n) - w(t_n)) = (e^n_w)^\ell + (R_h w(t_n) - w(t_n)), \\
			\partial^\tau_q (u^n_h)^\ell - \dot u(t_n)=&\  (\dot e^n_u)^\ell + (\partial^\tau_q - \d / \d t)R_h u(t_n)  + (R_h \dot u(t_n) - \dot u(t_n)).
		\end{aligned} 
	\end{equation}
	For the last terms we use the Ritz map error estimates of Lemma~\ref{Lemma:RitzMapError}, to estimate for $0 \leq j \leq q+1$,
	\begin{equation} \label{eq:MainProofRitzmaperror1}
		\begin{aligned}
			\| R_h u^{(j)} - u^{(j)} \|= \| R_h u^{(j)} - u^{(j)} \|_{H^1(\Om)} + \| \gamma (R_h u^{(j)} - u^{(j)})\|_{H^1(\Ga)} \leq ch, \\
			| R_h u^{(j)} - u^{(j)} |= \| R_h u^{(j)} - u^{(j)} \|_{L^2(\Om)} + \| \gamma (R_h u^{(j)} - u^{(j)})\|_{L^2(\Ga)} \leq ch^2 .
		\end{aligned}
	\end{equation}
	
	Therefore, with Lemma~\ref{Lemma:PeanoKernel}, we also obtain
	\begin{equation}
	\label{eq:MainProofRitzmaperror2}
		\begin{aligned}
			\| (\partial^\tau_q - \d / \d t)R_h u(t_n) \| \leq C \, \tau^q (1+h) , \andquad
			| (\partial^\tau_q - \d / \d t)R_h u(t_n) | \leq C \, \tau^q (1+h^2).
		\end{aligned}
	\end{equation}
	The same estimates hold for $w$. 
	
	For the first terms we use a combination of the stability part, Proposition~\ref{Prop:StabilityEst}, with the bounds for the defects in Proposition~\ref{Prop:ConsistencyEst} and the step size restriction, to arrive at
	\begin{align}
		\Big(\| e^n_u\|_h^2 + \| e^n_w \|_h^2 + \tau \sum_{k=q}^n \| \dot e^k_u \|_h^2\Big)^{\frac{1}{2}} \leq C \, (h^2+\tau^q). \nonumber
	\end{align}
	Then, we use the norm equivalence, see Remark~\ref{rem:normequivalence}, to obtain	
	\begin{align}
		\Big(\| (e^n_u)^\ell \|^2 + \| (e^n_w)^\ell \|^2 + \tau \sum_{k=q}^n \| (\dot e^k_u)^\ell \|^2\Big)^{\frac{1}{2}} \leq&\  C \, (h^2+\tau^q). \label{eq:MainProof2}
	\end{align}
	Combining the estimates \eqref{eq:MainProofdecomp}, \eqref{eq:MainProofRitzmaperror1}, \eqref{eq:MainProofRitzmaperror2}, and \eqref{eq:MainProof2}, finishes the proof. \hfill 
\end{proof}

\section{Numerical experiments}
\label{section:NumericalExperiments}

In this section we illustrate and complement our theoretical results with numerical simulations.

We will report on three numerical experiments:
\begin{itemize}
	\item A convergence experiment in the unit disk using a manufactured solution.
	\item A convergence experiment in the three-dimensional unit ball using a manufactured solution.
	\item We repeat some computations from \cite{KnopfLamLiuMetzger2021} with various initial data, allowing for a direct comparison.
\end{itemize}

\subsection{Convergence experiment in the unit disk}
\label{section:conv experiment - 2D}
We start by reporting on the convergence rate in time and space for a problem with a manufactured solution by adding inhomogeneities on the right-hand side, similar to the numerical experiments performed in \cite[Section~9]{HarderKovacs2022} for the spatial error.

Consider the Cahn--Hilliard equation with Cahn--Hilliard-type dynamic boundary conditions \eqref{eq:strong2CH} on the unit disk $\Om = \{x \in \R^2 \mid \|x\|_2 \leq 1 \}$, with parameters $\eps=\delta=\kappa=1$, and with the double-well potentials 
\begin{align*}
	W_\Om(u)=&\ \frac{1}{4}(u^2-1)^2, \qquad W_\Ga(u)= \frac{1}{4}(u^2-1)^2.	
\end{align*}

We use linear bulk--surface finite elements to discretize in space and the 3-step BDF method to discretize in time. In order report on the $L^2(\Om;\Ga)$- and $H^1(\Om;\Ga)$-norm errors, we construct  inhomogeneities for each equation, such that the exact solutions are known to be
\begin{align*}
	u(t,x_1,x_2)= \exp(-t) x_1x_2, \andquad w(t,x_1,x_2)= \exp(-t) x_1x_2.
\end{align*}

The modified system is then given by: Find $u,w \colon \bar \Om \times [0,T] \to \mathbb{R}$ satisfying
\begin{subequations} \label{eq:Modstrong2CH}
	\begin{alignat}{2}
		\dot u =&\ \varDelta w + f_\Om&&\text{ in } \Om, \\
		w=&\ - \eps \varDelta u + \frac1\eps W'_\Om (u) +g_\Om&&\text{ in } \Om,\\
		\dot u =&\ \varDelta_\Ga w - \partial_\nu w +f_\Ga&&\text{ on } \Ga, \\
		w=&\ - \delta \kappa \varDelta_\Ga u + \frac1\delta W'_\Ga (u) + \eps \partial_\nu u +g_\Ga\qquad &&\text{ on } \Ga,
	\end{alignat}
\end{subequations}
where for $\varphi(t,x_1,x_2)=\exp(-t) x_1x_2$ the inhomogeneities are defined as
\begin{align}
	f_\Om:=-\varphi, \quad g_\Om:=\varphi-\frac1\eps W_\Om'(\varphi), \quad f_\Ga:= 5\varphi, \quad g_\Ga:=(1-4\delta\kappa -2\epsilon)\varphi - \frac1\delta W_\Ga'(\varphi).
\end{align}

\begin{rem}
	By a minimal extension the problem \eqref{eq:Modstrong2CH}, fits into the framework of the theoretical result Theorem~\ref{Thm:ConvergenceMainResult}. The error equations are not changed by this modification, only four simple defect terms are added and we need $L^2(\Om)$-regularity in the bulk, for $f_\Om$ and $g_\Om$, and $L^2(\Ga)$-regularity on the surface, for $f_\Ga$ and $g_\Ga$, to obtain the analogous regularity result.
\end{rem}

For given meshes with degrees of freedom $2^k \cdot 20$, for $k \in \{0,\dotsc,12\}$, and given time step sizes $\tau \in \{0.25, 0.125, 0.05, 0.025, 0.0125, 0.005, 0.0025, 0.00125\}$, the starting values are taken to be the nodal interpolation of the exact solution. 
\begin{figure}[!t]
	\centering\includegraphics[width=\textwidth]{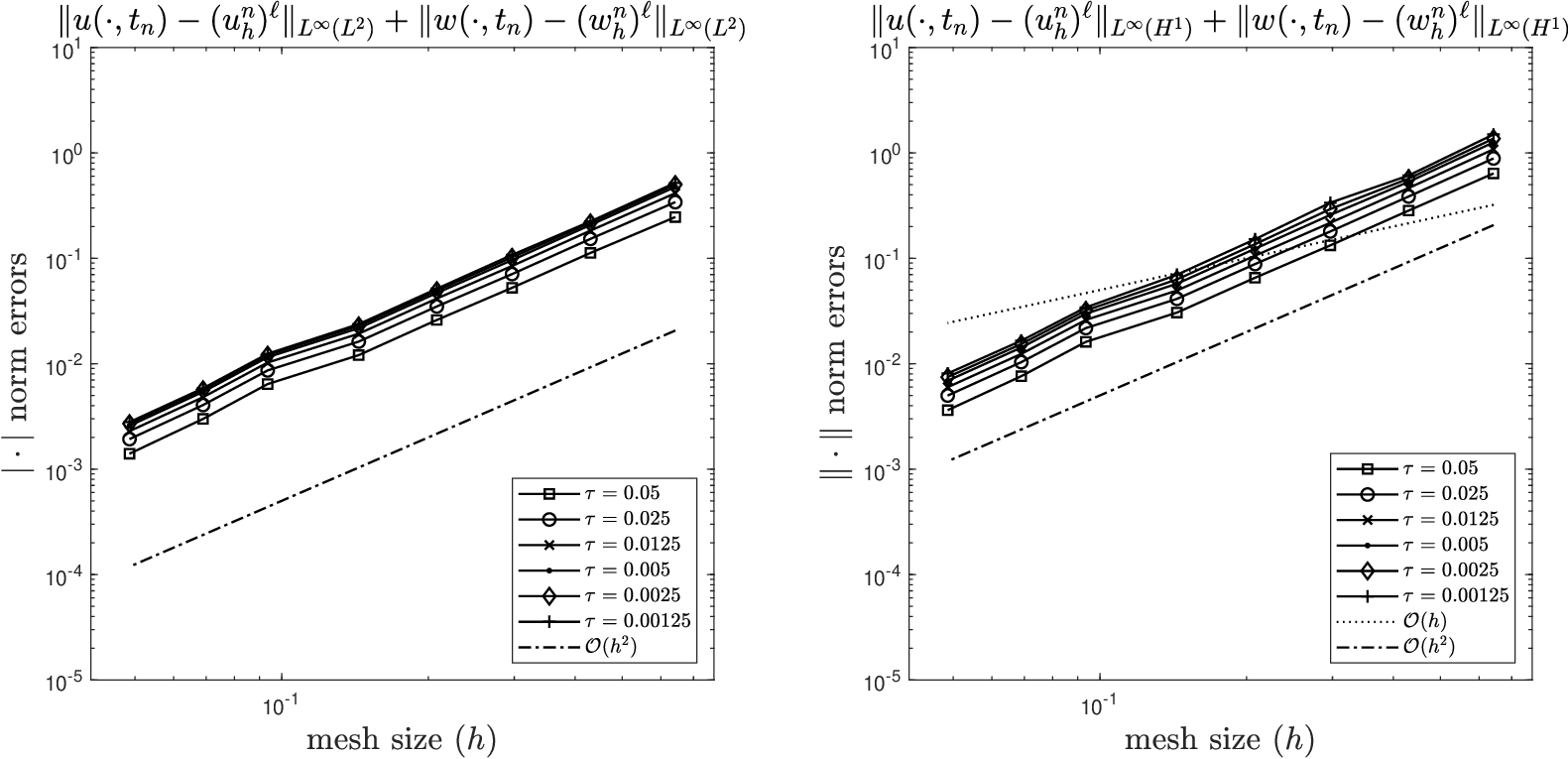}
	\caption{Spatial convergence plots for the linear bulk--surface FEM/BDF3 approximation to the Cahn--Hilliard equation with Cahn--Hilliard-type dynamic boundary conditions in the two-dimensional unit disk with free energy double well potentials.}
	\label{fig:SpatialConvergence}
\end{figure}

In Figure~\ref{fig:SpatialConvergence} we report on the $L^\infty([0,1],L^2(\Om;\Ga))$ norm (left) and the $L^\infty([0,1],H^1(\Om;\Ga))$ norm (right) of the errors between the numerical solution and the interpolation of the phase-field $u$ and the chemical potential $w$. 
The logarithmic plots show the errors against the mesh width $h$, for different time step sizes $\tau$. We observe, that the spatial error dominates for small time step sizes, and is of order $h^2$ for the $L^\infty([0,1],L^2(\Om;\Ga))$ norm, as proved in Theorem~\ref{Thm:ConvergenceMainResult}. For the $L^\infty([0,1],H^1(\Om;\Ga))$ norm, we even observe a better convergence rate of $h^2$, than the linear convergence rate obtained in Theorem~\ref{Thm:ConvergenceMainResult} (due to the super-convergence when comparing with the interpolation of the exact solutions). Note that similar convergence test can be found in \cite[Section~9]{HarderKovacs2022}.

\begin{figure}[!t]
	\centering\includegraphics[width=\textwidth]{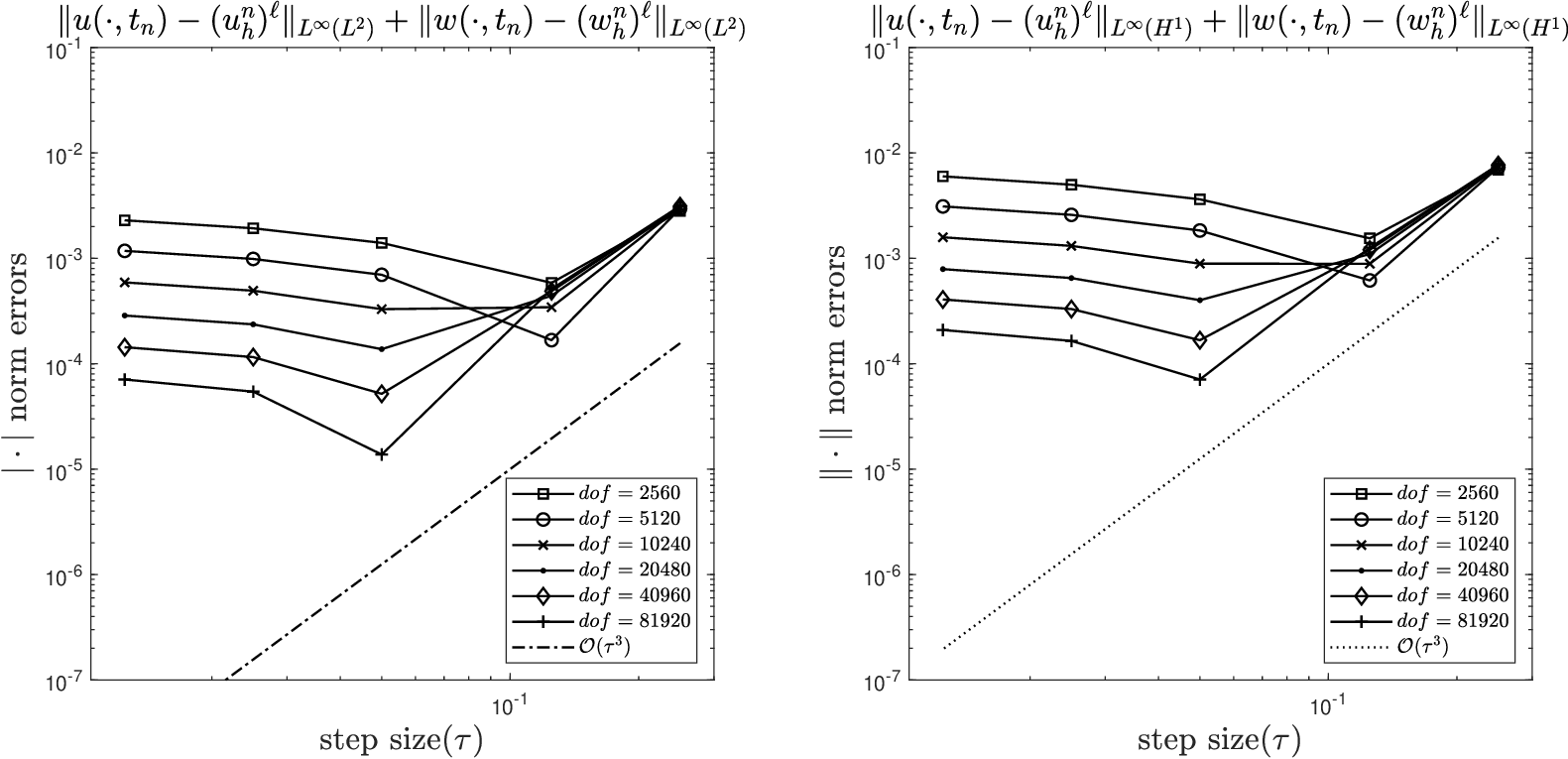}
	\caption{Temporal convergence plots for the linear bulk--surface FEM/BDF3 approximation to the Cahn--Hilliard equation with Cahn--Hilliard-type dynamic boundary conditions in the two-dimensional unit disk with free energy double well potentials.}
	\label{fig:TemporalConvergence}
\end{figure}

In Figure~\ref{fig:TemporalConvergence} the same errors are displayed, now against the time step size $\tau$, and with different mesh sizes. We again observe that, for small enough mesh size $h$, the temporal error dominates and admits the convergence rate $\tau^3$ for the $3$-step BDF method, as predicted by our Theorem~\ref{Thm:ConvergenceMainResult}.

\subsection{Convergence experiment in the unit sphere}
\label{section:conv experiment - 3D}

We also performed a convergence experiment in the three-dimensional unit ball, analogous to the one in Section~\ref{section:conv experiment - 2D}. 
The setup of the experiment is identical, except the inhomogeneities are now constructed, such that the exact solutions are known to be
\begin{align*}
	u(t,x)= \exp(-t) x_1 x_2 x_3, \andquad w(t,x)= \exp(-t) x_1 x_2 x_3 .
\end{align*}

\begin{figure}[!t]
	\centering\includegraphics[width=\textwidth]{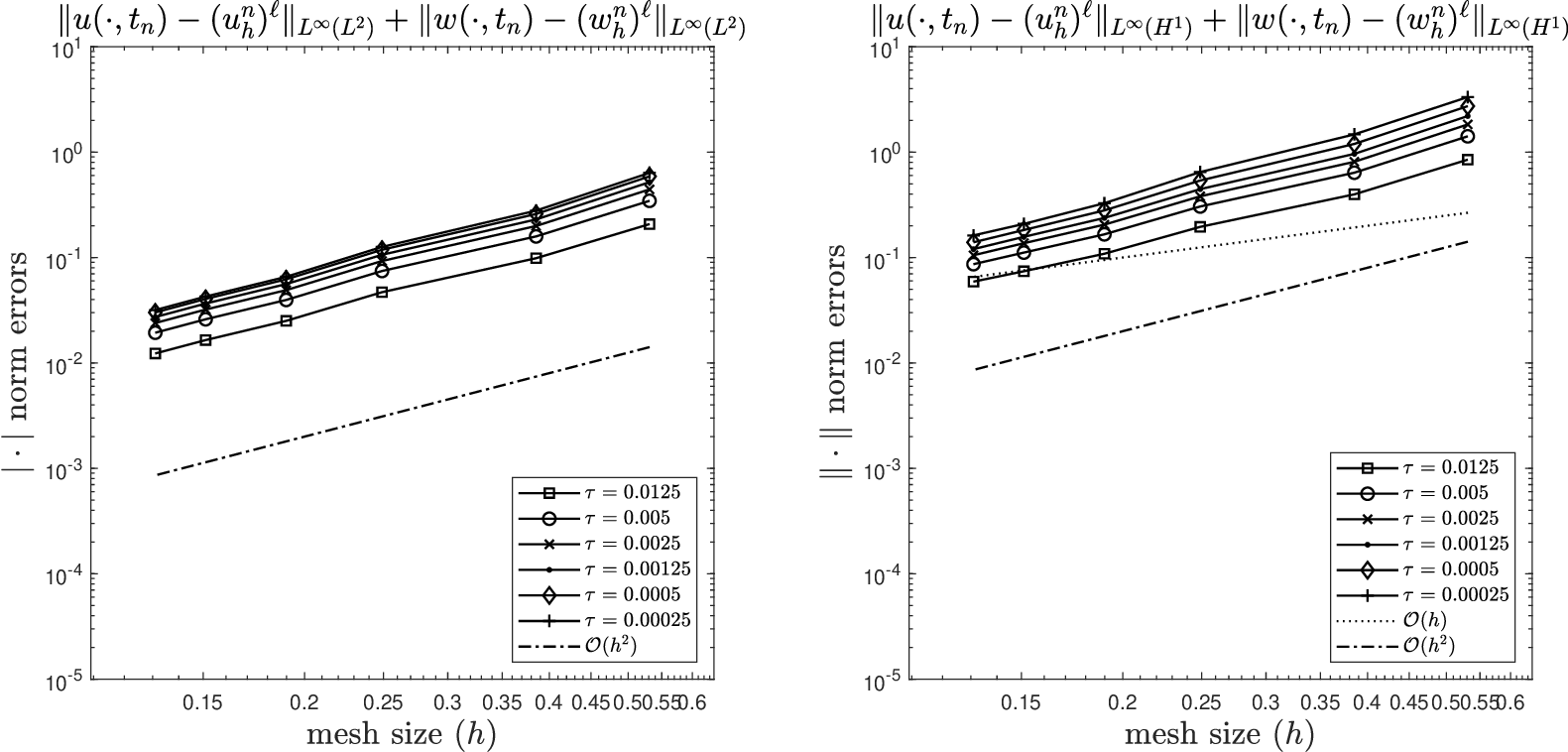}
	\caption{Spatial convergence plots for the linear bulk--surface FEM/BDF2 approximation to the Cahn--Hilliard equation with Cahn--Hilliard-type dynamic boundary conditions in the three-dimensional unit ball with free energy double well potentials.}
	\label{fig:SpatialConvergence - 3D}
\end{figure}

In Figure~\ref{fig:SpatialConvergence - 3D} we report on the errors for the three-dimensional test-problem in the $L^\infty([0,1],L^2(\Om;\Ga))$ norm (left) and the $L^\infty([0,1],H^1(\Om;\Ga))$ norm (right), for meshes with degrees of freedom $172,587,1766,4385,8480,13574$ and time step sizes $\tau \in \{0.0125, 0.005, 0.0025, 0.00125, 0.0005, 0.00025\}$. We have observed that in this three-dimensional problem the spatial error is heavily dominating the temporal error, so we would need extremely fine meshes to observe the temporal convergence order, hence we only report on spatial convergence.

\begin{figure}[!t]
	\centering\includegraphics[width=\textwidth,clip]{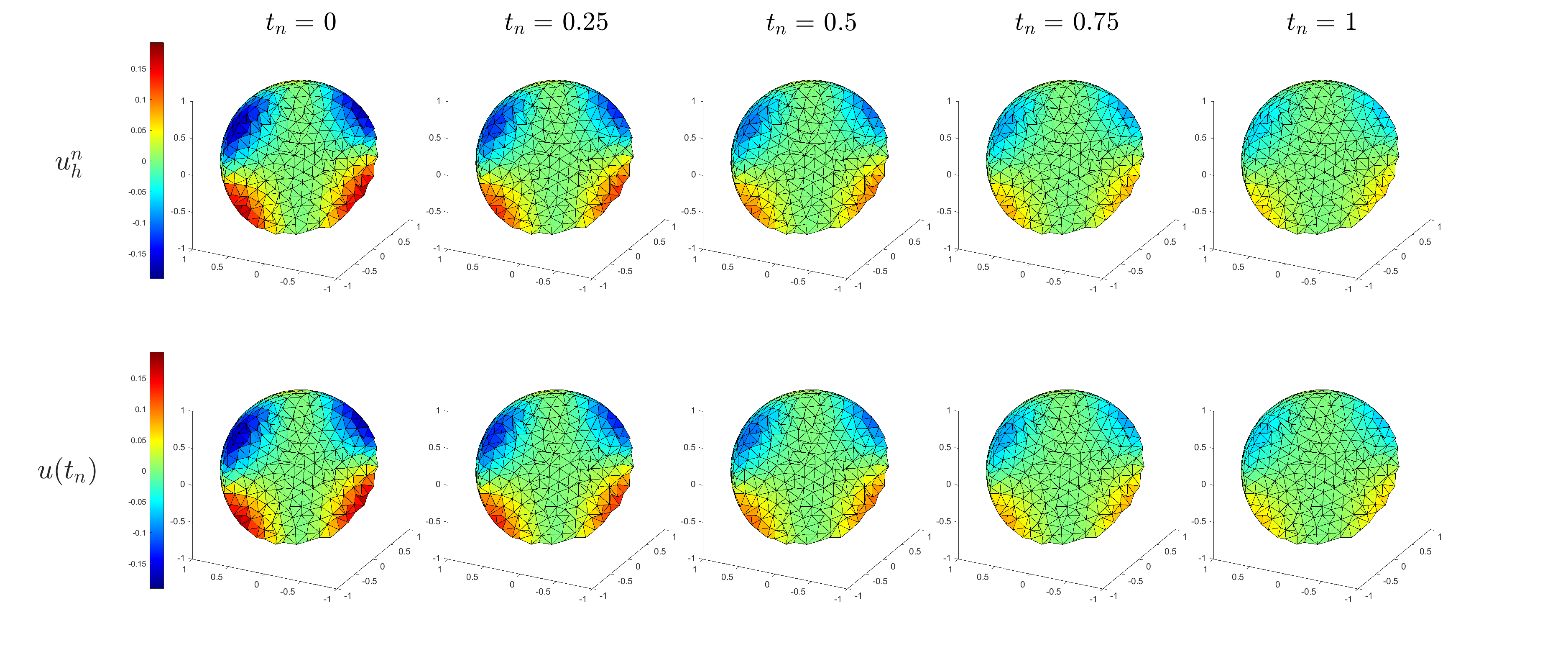}
	\caption{Evolution plot on a cross-section of the three-dimensional unit ball, for the exact solution and numerical approximation by the bulk-surface finite element/BDF2 method of the convergence experiment, with $\tau=0.001$ and a mesh with dof$=1766$.}
	\label{fig:EvolutionPlot - 3D}
\end{figure}

Additionally, in Figure~\ref{fig:EvolutionPlot - 3D} we report on the time evolution of the numerical approximation to $u$ at various times. The plots display a vertical cross-section of the domain.

\subsection{Numerical experiments on the evolution of three different initial data}
\begin{figure}[!t]
	\centering\includegraphics[width=\textwidth]{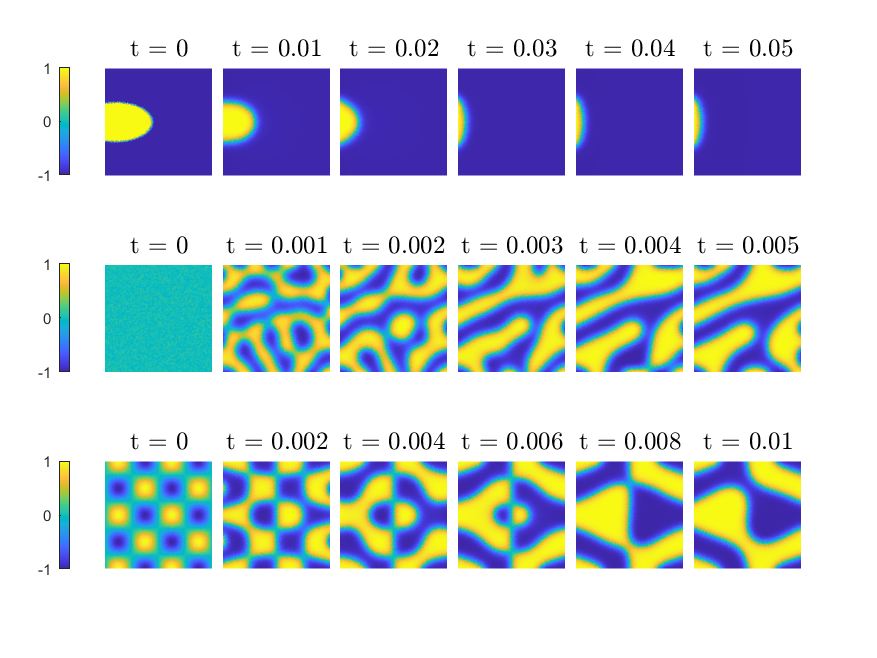}
	\caption{Evolution plot for three different kind of initial data by the bulk-surface finite element/BDF2 method.}
	\label{fig:EvolutionPlot}
\end{figure}

A simulation of the phase separation of the Cahn--Hilliard equation with dynamic boundary conditions with diffusion interface parameters $\epsilon=\delta =0.02$ and $\kappa=1$, for three different initial data is presented in Figure~\ref{fig:EvolutionPlot}. We again set the non-linear double-well potentials
\begin{align*}
	W_\Omega(u)=&\ \frac{1}{8}(u^2-1)^2, \qquad W_\Gamma(u)= \frac{1}{8}(u^2-1)^2.
\end{align*} 
These choices allow a direct comparison with the numerical experiments of \cite{KnopfLamLiuMetzger2021}, wherein the authors use a semi-implicit Euler method, with convex-concave splitting for the nonlinear potential.

We approximate the solution by our numerical scheme, combining the linear finite element method and the linearly implicit BDF method of order $2$, with a triangular mesh consisting of $11660$ nodes and $\tau=10^{-5}$. The three different initial data are given on the domain $\Omega=(0,1)^2 \subset \mathbb{R}^2$ by:
\begin{enumerate}
	\item An elliptical shaped droplet centred at $(0.1,0.5)$, a maximal horizontal elongation of $0.6814$, and a maximal vertical elongation of $0.367$, see \cite[Section 6]{KnopfLamLiuMetzger2021}.
	\item A uniformly distributed random value in $[-0.1,0.1]$ for any node, see \cite[Section 7.2.2]{KnopfLam2020}.
	\item The initial data $u^0(x)=\sin(4\pi x_1) \cos(4\pi x_2)$, see \cite[Section 7.2.2]{KnopfLam2020}.
\end{enumerate}
The starting value $\bfu^1$ is computed by a linearly implicit Euler step (which has local order $2$ in time). Each row of Figure~\ref{fig:EvolutionPlot} displays the evolution of phase-field $u$, depicted in the first column, plotted at different time-steps.

\begin{figure}[!t]
	\centering\includegraphics[width=\textwidth]{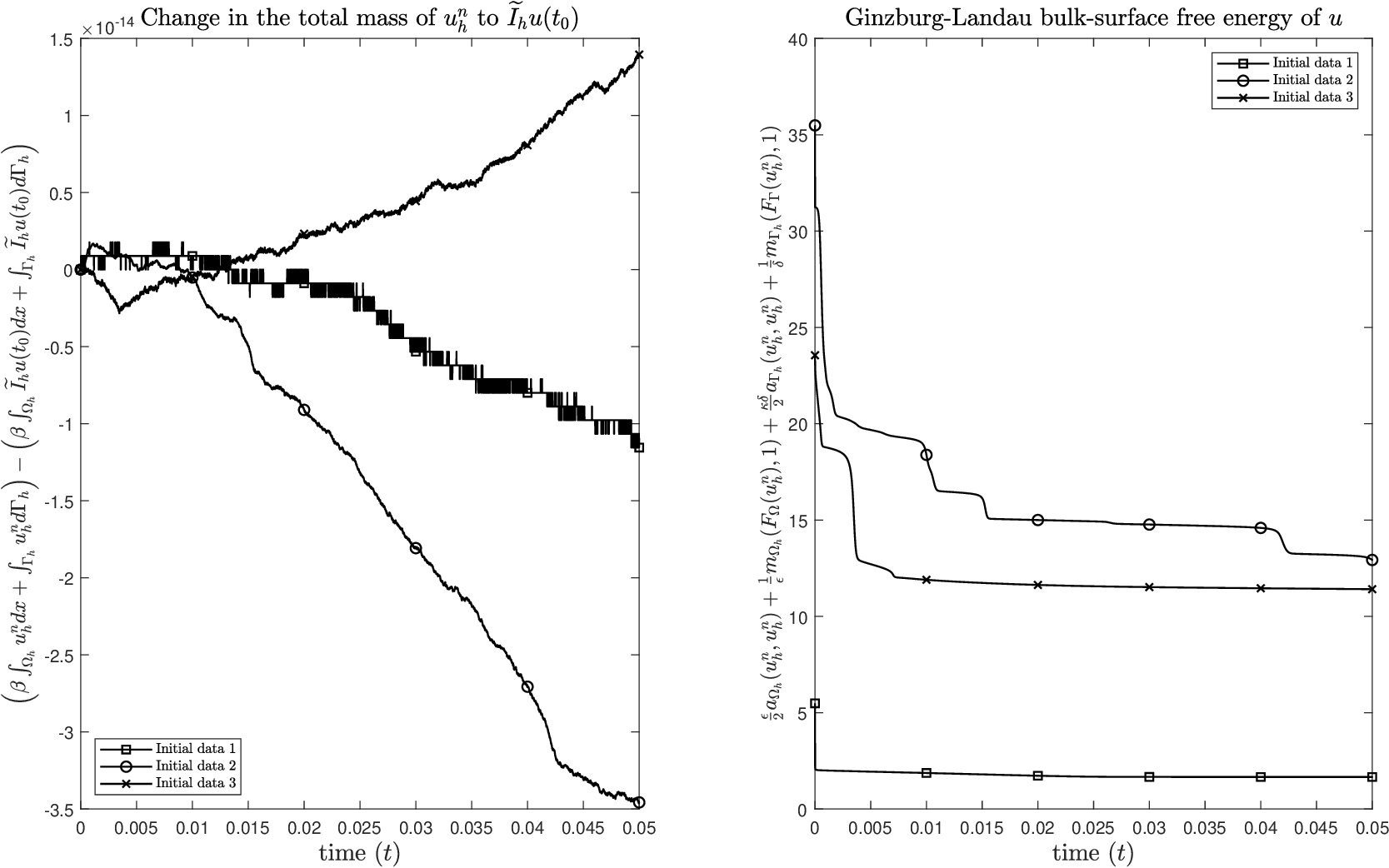}
	\caption{Mass and energy plot for three different kind of initial data by the bulk-surface finite element/BDF2 method (note axis-scaling).}
	\label{fig:MassEnergyPlot}
\end{figure}

In Figure~\ref{fig:MassEnergyPlot} we report on the evolution of the discrete total mass and discrete Ginzburg-Landau bulk-surface free energy, see \eqref{eq:GinzLand}, for the three different initial data. We observe a mass conservation, as discussed in Remark~\ref{rem:MassEnergy}, and a dissipation of the energy.

\section*{Acknowledgments}
The work of Bal\'azs Kov\'acs is funded by the Heisenberg Programme of the Deutsche Forschungsgemeinschaft (DFG, German Research Foundation) -- Project-ID 446431602,
and by the DFG Research Unit FOR 3013 \textit{Vector- and tensor-valued surface PDEs} (BA2268/6–1).

\appendix
\section{$H^2$-regularity via \cite{ElliottRanner2015}}
\label{appendix}

We will now present a compact proof to the $H^2$-regularity result. For the analogous version for singular potentials we refer to \cite[Theorem~2.2, Section~4.4]{ColliFukao_2015}.

\begin{thm}
\label{theorem:H^2 regularity - appendix}
	Let $u_0 \in H^2(\Om) \cap H^2(\Ga)$ and $(u, w)$ be the solution pair of \eqref{eq:weakCH}. Then $u \in L^\infty([0,T];H^2(\Om;\Ga))$ and $w \in L^2([0,T];H^2(\Om;\Ga))$, with the bounds
	\begin{equation}
		\sup_{t \in (0, T)} \|u\|^2_{H^2(\Om;\Ga)} + \int_0^T \|w\|^2_{H^2(\Om;\Ga)} \,dt \leq C_2(u_0).
	\end{equation}
\end{thm}

\begin{proof}
(a) The first step is to prove the improved estimates analogous to \cite[Lemma~4.2]{ElliottRanner2015} for the semi-discrete solution $(u_h,w_h)$. 
First, by standard ODE theory the short-time well-posedness of the semi-discrete solution $(u_h,w_h)$ follows and also that both functions are $C^1$ in time, and by an energy argument these can be extended to $[0,T]$, see \cite[Lemma~4.2]{ElliottRanner2015}. Second, by using energy estimates testing \eqref{eq:sdcCHw1} by $\dot u_h$ and the time derivative of \eqref{eq:sdcCHw2} by $w_h$, similar to Part (ii) above, the analogue of \cite[Lemma~4.2]{ElliottRanner2015} follows.

From this improved estimates, there follows the following energy-type estimate 
\begin{equation}
	\int_0^T \|\partial_t u\|^2_{L^2(\Om;\Ga)} \d t + \sup_{t \in (0, T)} \|w\|^2_{L^2(\Om;\Ga)} \leq C_2(u_0) .
\end{equation}

(b) Using a version of the Aubin--Nietsche trick, we translate the fact that $(u, w)$ are solutions of \eqref{eq:absCH} into
\begin{equation*}
	a(u, \varphi^w) = m(f_1, \varphi^w) \andquad a(w, \varphi^u) = m(f_2, \varphi^u) \qquad \text{ for all } \varphi^u, \varphi^w \in V ,
\end{equation*}
where 
\begin{equation*}
	f_1 := w - W'(u) \andquad f_2 = \dot u .
\end{equation*}

The above improved bounds, combined with the standard energy bounds, give that  $f_1 \in L^\infty([0,T];L^2(\Om;\Ga))$ and $f_2 \in L^2([0,T];L^2(\Om;\Ga))$. Standard theory of elliptic partial differential equations, e.g.~\cite{Aubin_book}, gives $u \in L^\infty([0,T];H^2(\Om;\Ga))$ and $w \in L^2([0,T];H^2(\Om;\Ga))$. The estimate in the theorem is shown by combining previous bounds on $f_1$ and $f_2$.
\hfill 
\end{proof}

\bibliographystyle{abbrvnat}
\bibliography{C-Heq_review}

\end{document}